\documentclass{article}
\usepackage{graphicx}
\usepackage{amsfonts}
\usepackage{floatrow}
\usepackage{authblk}
\usepackage{cite}
\usepackage{bm}
\usepackage{amsmath}
\usepackage{amssymb}
\renewcommand\subfloatrowsep{\hskip -5\columnsep}
\newcommand{\mat}[2]{\left[ \begin{array}{*{#1}{c}}#2\end{array}\right]}

\begin{document}
\thispagestyle{empty}
\vspace*{-0.5in}
\begin{center} {\large{ \bf Semi-invariants of binary forms and symmetrized graph-monomials}}
\end{center}
\begin{center} {\bf Shashikant Mulay  \\ \small Department of Mathematics,  
University of Tennessee \\ Knoxville, TN 37996 U. S. A.;  e-mail: smulay@utk.edu \\ May 17, 2019} \end{center}

\vspace{0.2in}

\noindent \underline{ \bf Abstract}: This article provides a method for constructing invariants and 
semi-invariants of a binary $N$-ic form over a field $k$ characteristics $0$ or $p > N$. A practical 
and broadly applicable sufficient condition for ensuring nontriviality of the symmetrization of a 
graph-monomial is established. This allows construction of infinite families of invariants (especially, 
skew-invariants) and families of $k$-linearly independent semi-invariants. These constructions are very 
useful in the quantum physics of Fermions. Additionally, they permit us to establish a new polynomial-type 
lower bound on the coefficient of $q^{w}$ in $(q - 1) {N + d \choose d}_{q}$ for all sufficiently large 
integers $d$ and $w \leq N d / 2$.  \\



\vspace{.2in}

Fix an integer $N \geq 2$. Let $k$ be a field of characteristic either $0$ or strictly greater than $N$. 
Let $X$, $Y$, $t$, $z_{1}, \ldots , z_{N}$ be indeterminates. Let $E_{1}(t), \ldots , E_{N}(t)$ and $f(X + t)$ 
be the polynomials defined by
\[f(X + t) \; := \; \prod _{i=1}^{N} (X + z_{i} + t) \; =: \; X^{N} + \sum _{i=1}^{N} E_{i} (t) X^{N - i} .\]
For $1 \leq i \leq N$, let $e_{i} := E_{i} (0)$. Then, $f(X) = X^{N} + e_{1} X^{N - 1} + \cdots + e_{N}$. 
A polynomial $P(e_{1}, \ldots , e_{N}) \in k[e_{1}, \ldots , e_{N}]$ is said to be {\it translation invariant}
provided $P(E_{1} (t), \ldots , E_{N} (t)) = P(e_{1}, \ldots , e_{N})$. It is a (well known) simple exercise to
verify that the subring $k[y_{1}, \ldots , y_{N - 1}]$ of $k[e_{1}, \ldots , e_{N}]$, where $y_{i} := E_{i} (-e_{1} / N)$ 
for $1 \leq i \leq N$, is the ring of all translation invariant members of $k[e_{1}, \ldots , e_{N}]$. Furthermore,
we have $k[y_{1}, \ldots , y_{N - 1}] = k[e_{1}, \ldots , e_{N}] \cap k[z_{1} - z_{2}, \ldots , z_{1} - z_{N}]$ 
({\it e.g.}, see Ch. 2, Theorem 1 of [11]). 
A polynomial $h \in k[e_{1}, \ldots , e_{N}]$ is said to be homogeneous of {\it weight $w$} provided as a polynomial in 
$z_{1}, \ldots , z_{N}$, $h$ is homogeneous of degree $w$. Note that $y_{i}$ is homogeneous of weight $i + 1$ for
$1 \leq i \leq N$. Next, consider the (generic) binary form $F := \sum a_{i} X^{i} Y^{N - i}$ of degree $N$ where 
$a_{0}$ is an indeterminate and $a_{i} := a_{0} e_{i}$ for $1 \leq i \leq N$. A {\it semi-invariant of $F$ of degree $d$
and weight $w$} is a polynomial $Q \in k[a_{0}, a_{1}, \ldots , a_{N}]$ such that $Q = a_{0}^{d} P(e_{1}, \ldots , e_{N})$
where $P(e_{1}, \ldots , e_{N})$ is translation invariant, homogeneous of weight $w$ and has total degree $\leq d$ in 
$e_{1}, \ldots , e_{N}$. For $0 \leq i \leq N$, the weight of $a_{i}$ is defined to be $i$. Then, note that $Q$ is 
homogeneous of degree $d$ and weight $w$ in $a_{0}, \ldots , a_{N}$. An {\it invariant} of $F$ of degree $d$ is a 
semi-invariant of $F$ of degree $d$ and weight $N d / 2$. For a fixed $N$, the set of semi-invariants (of
the binary $N$-ic $F$) of degree $d$ and weight $w$ form a finite dimensional $k$-linear subspace of 
$k[a_{0}, a_{1}, \ldots , a_{N}]$. This subspace is known to be trivial unless $2 w \leq N d$. Provided $\mbox{char} \,k = 0$  
and $2 w \leq N d$, a theorem of Cayley-Sylvester proves that the dimension of the aforementioned space 
of semi-invariants of degree $d$ and weight $w$ is the coefficient of $q^{w}$ in $(q - 1) {N + d \choose d}_{q}$ where 
${N + d \choose d}_{q}$ is the $q$-binomial coefficient (see [6], [18] or Theorem 5 of [11]). Let $p_{w}(N, d)$ denote
the coefficient of $q^{w}$ in $q^{w}$ in ${N + d \choose d}_{q}$. Then, $p_{w}(N, d)$ is the number of integer-partitions
of $w$ in at most $N$ parts with each part $\leq d$. As a corollary of the Cayeley-Sylvester theorem, we then have
$p_{w}(N, d) \geq  p_{w - 1}(N, d)$ for $2 \leq w \leq N d / 2$; this establishes {\it unimodality} of the coefficients of 
${N + d \choose d}_{q}$. Since $p_{w}(N, d) - p_{w - 1}(N, d)$ are the dimensions of spaces of semi-invariants, it is natural
to investigate explicit (lower, upper) bounds on them. Recently, some interesting lower bounds on $p_{w}(N, d) - p_{w - 1}(N, d)$
have come to light (see [4], [12], [19] and their references). This article has two objectives: provide explicit methods of
constructing a class of $k$-linearly independent semi-invariants and obtain a new lower bound on $p_{w}(N, d) - p_{w - 1}(N, d)$
for certain pairs $(w, d)$. The non-trivial lower bounds of [4], [12] and [19] are valid for $\min \{N, d \} \geq 8$ but
for all sufficiently large values of $d$ and $w$, they do not depend on $(w, d)$. In contrast, our lower bounds (see Theorem 3) 
are polynomials in $w$ for all $(N, d)$; Example 3, 4 and Remark 5 appearing at the end of the article present a more detailed comparison. 
In this article, we investigate the algebra of semi-invariants; not the combinatorics of $q$-binomial coefficients. In the rest of 
the introduction, we describe our motivation for, and our method of, constructing semi-invariants of a binary $N$-ic form. 

Ever since the theory of invariants of binary forms was founded, invariant-theorists have explored and devised methods for 
writing down concrete invariants; however, each of these methods has its own shortcomings. The `symbolic method' of classical 
invariant theory (see [3], [6], [7], [9]) provides an easy recipe for formulating symbolic expressions that yield invariants 
and semi-invariants. But, without full expansion (or un-symbolization) one does not know whether a given symbolic expression 
yields a {\em nonzero} semi-invariant. Here we prefer the other method, {\it i.e.}, the method of symmetrized graph-monomials. 
This too was known to classical invariant theorists (see [13], [14], [17]). It poses the problem of finding a useful criterion 
to determine nonzero-ness of the symmetrization. Historically, Sylvester and Petersen considered this problem; in fact, Petersen 
formulated a sufficient (but not necessary) condition that ensures zero-ness of the symmetrization. For a detailed historical 
sketch of this topic, we refer the reader to [16]. In [16], nonzero-ness of the symmetrization of a graph-monomial is shown to be 
equivalent to certain properties of the orientations and the orientation preserving graph-automorphisms of the underlying graph; 
but as matters stand, verification of these properties is as forbidding as is a brute force computation of the desired symmetrization. 
Our interest in {\em construction}, as opposed to {\em existence}, of invariants and semi-invariants stems primarily from the need to 
obtain explicitly described {\it trial wave functions} for systems of $N$ strongly correlated Fermions in fractional quantum Hall state. 
Such a trial wave function is essentially determined by a so called {\it correlation function}. The intuitive approach of physics 
presents such a correlation function as a symmetrization of a {\it monomial} obtained from the graph of correlations representing 
allowed strong interactions between $N$ Fermions. It so happens that this correlation function turns out to be a semi-invariant 
(an invariant in certain cases), of a binary $N$-ic form. In this article, we establish an easy to use yet broadly applicable 
sufficient criterion (see Theorem 1) for non-triviality of a symmetrized graph-monomial. Besides enabling explicit 
constructions of the desired trial wave functions, Theorem 1 is also interesting from a purely invariant theoretic point of view. 
Following Theorem 1, we exhibit a sample of its applications (see Theorem 2, Theorem 3).. 
 
A {\it multigraph} is a graph in which multiple edges are allowed between the same two vertices of the graph.
Consider a loopless undirected multigraph $\Gamma$ on finitely many (at least two) vertices labeled 
$1, 2, \ldots , N $; multigraph $\Gamma$ is said to be {\it $d$-regular} provided each vertex of $\Gamma$ 
has the same degree $d$. In the figures below, $\Gamma _1$ is seen to be a $2$-regular multigraph and 
the multigraphs $\Gamma _2$, $\Gamma _3$ both are $3$-regular. 

\begin{figure}[h!]
\hspace*{-2.5cm}
\begin{floatrow}
 \ffigbox{\includegraphics[scale = 0.25]{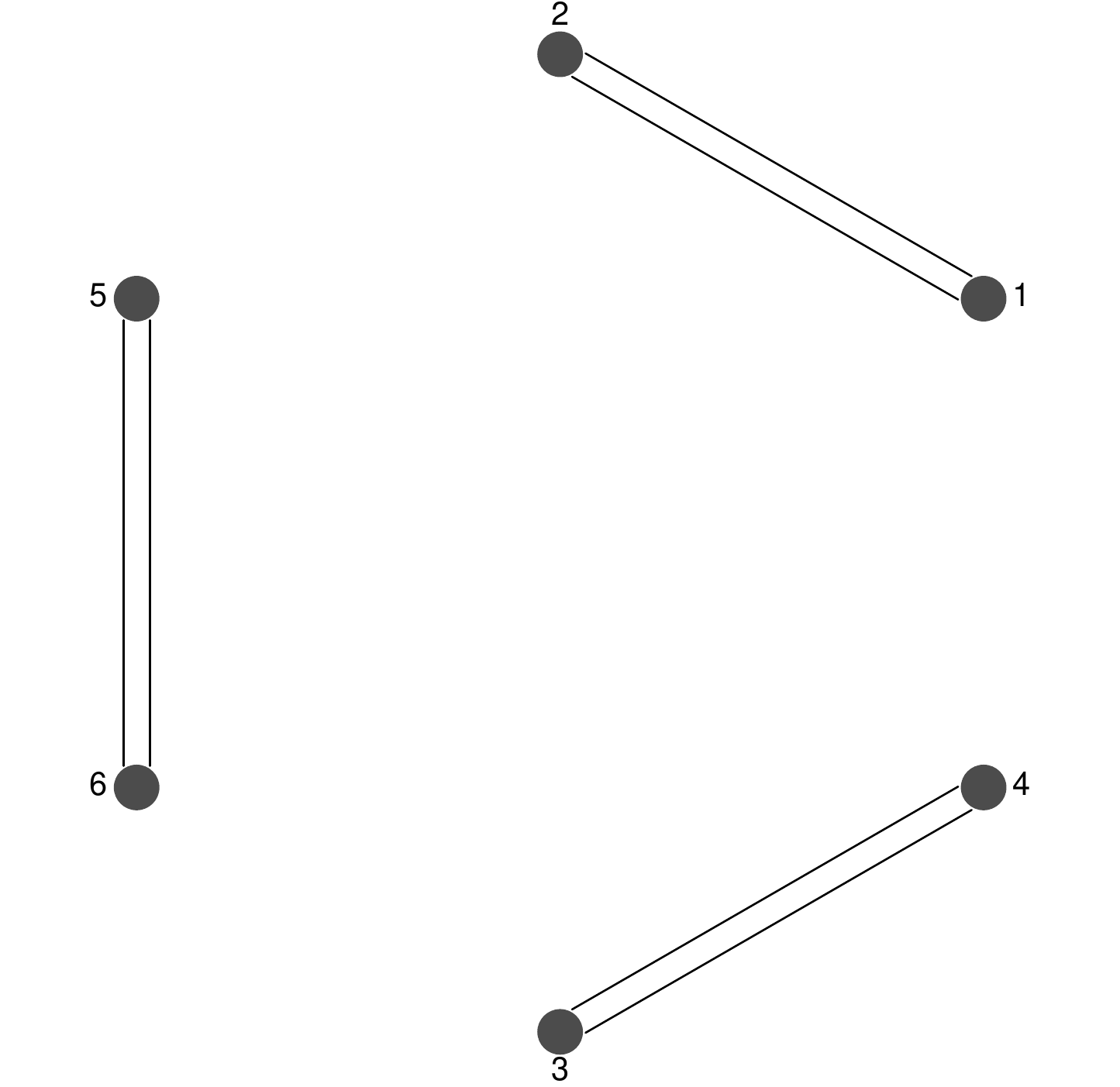}}{\caption{$\Gamma _1$}} 
 \subfloatrowsep
 \ffigbox{\includegraphics[scale = 0.25]{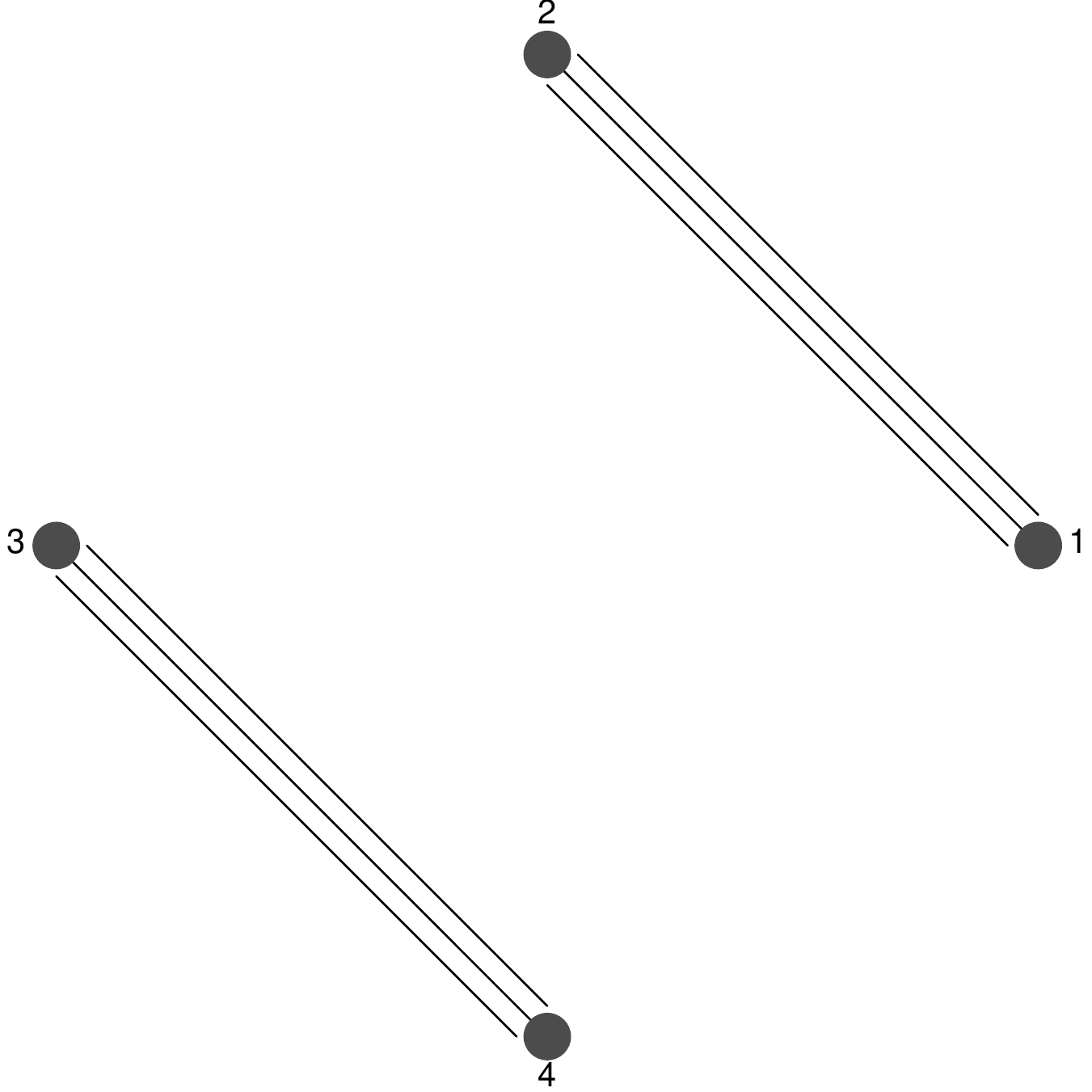}}{\caption{$\Gamma  _2$}} 
 \subfloatrowsep
 \ffigbox{\includegraphics[scale = 0.25]{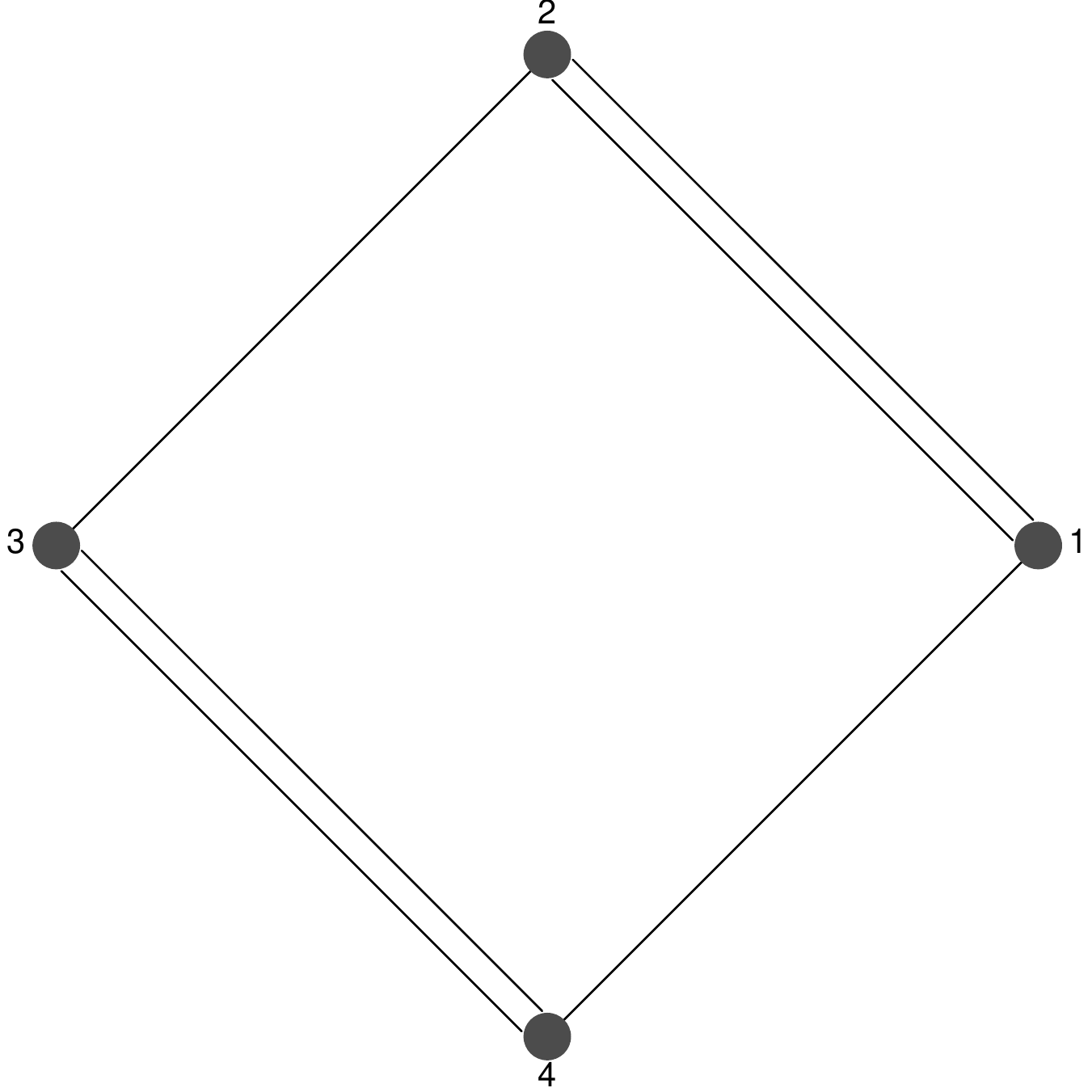}}{\caption{$\Gamma  _3$}} 
\end{floatrow}
\end{figure}

Let $\varepsilon ( \Gamma, i, j)$ be the number of edges in $\Gamma$ connecting vertex $i$ to vertex $j$.
The {\it graph-monomial} of $\Gamma$, denoted by $\mu (\Gamma)$, is the polynomial in $z_{1}, \ldots , z_{N}$ defined by 
\[ \mu (\Gamma) \; := \; \prod _{1 \leq i < j \leq N} (z_{i} - z_{j})^{\varepsilon (\Gamma, i, j)} . \]
Let $g(\Gamma)$ denote the {\it symmetrization} of $\mu (\Gamma)$, {\it i. e.}, $g(\Gamma) := \sum \mu_{\sigma} (\Gamma)$, 
where the sum ranges over the permutations $\sigma$ of $\{1, 2, \ldots , N \}$ and $\mu _{\sigma} (\Gamma)$ stands for the 
product of $(z_{\sigma (i)} - z_{\sigma (j)})^{\varepsilon (\Gamma, i, j)}$; $1 \leq i < j \leq N$. In the classical invariant 
theory of binary forms (where $k = {\Bbb C}$), it is well known that if $\Gamma$ is $d$-regular on $N$ vertices, then $g(\Gamma)$ 
is a (relative) invariant of degree $d$ (and weight $N d /2$) of the binary $N$-ic form $F$. Moreover, the 
vector space of invariants of $F$ of degree $d$ is spanned by the set of symmetrized graph monomials corresponding to the 
$d$-regular multigraphs on $N$ vertices (for a proof see [6] or its modern treatment: Ch. 2, Theorem 4 of [11]). If $\Gamma$ is 
not $d$-regular for any $d$, then $g(\Gamma)$ is a semi-invariant (as defined in [6], [7]) of $F$ irrespective of the characteristic 
of $k$. For example, $g(\Gamma _{1})$ is a quadratic invariant of a binary sextic (investigated in [5]) and each of $g(\Gamma _{2})$, 
$g(\Gamma _{3})$ is a cubic invariant of a binary quartic. It can be easily verified that $g(\Gamma _{2})$ is identically $0$ whereas 
$g(\Gamma _{3})$ is essentially the only nonzero cubic invariant of a binary quartic. In general, given a nonzero semi-invariant of $F$, 
there is no known method to determine whether the invariant is $g(\Gamma)$ for somme multigraph $\Gamma$. Also, for non-isomorphic 
multigraphs $\Gamma$, $\Gamma ^{\prime}$, their corresponding semi-invariants $g(\Gamma)$, $g(\Gamma ^{\prime})$ may be numerical 
multiples of each other. Clearly, it is desirable to understand the types of multigraph $\Gamma$ for which $g(\Gamma)$ is nonzero. 
For then, we get a natural method of constructing nonzero semi-invariants of $F$. 

In the physics of Fermion-correlations, vertices of $\Gamma$ correspond to Fermions and the edges in $\Gamma$ represent 
correlations (a repulsive interaction) between the Fermions; here, it suffices to work over ${\Bbb C}$. A multigraph 
$\Gamma$ is called a {\it configuration} of Fermions provided $g(\Gamma)$ is nonzero, and then $g (\Gamma)$ is called the 
correlation-function of this configuration. A configuration $\Gamma$ need not be $d$-regular for any $d$. In physics a 
configuartion $\Gamma$ is as important as its associated correlation function $g(\Gamma)$. This leads to some 
interesting new problems that do not seem to have any parallels in the theory of invariants. For example, let 
$p (\Gamma)$ and $L (\Gamma)$ denote the maximum of and the sum of all $\varepsilon (\Gamma, i, j)$ respectively. For 
fixed integers $N$, $L$ and $d$, consider the set $C (N, L, d)$ of multigraphs $\Gamma$ with the maximum vertex-degree 
$d$, $L (\Gamma) = L$ and $g (\Gamma) \neq 0$. Let $p(N, L, d)$ denote the minimum of $p (\Gamma)$ as $\Gamma$ ranges over 
$C(N, L, d)$. A configuration $\Gamma \in C (N, L, d)$ is {\it minimal} if $p (\Gamma) = p(N, L, d)$. It is known 
(see [11], [15]) that the lowest energy configurations (or states) $\Gamma$ are those with the least $p(\Gamma)$. Thus
one needs to estimate $p(N, L, d)$ for a given triple $(N, L, D)$. Likewise, given $\Gamma, \, \Gamma ^{\prime} \in C (N, L, d)$, it 
is of interest to know when $g (\Gamma)$ is (or is not) a constant multiple of $g(\Gamma ^{\prime})$. Without digressing into 
deeper physics, we simply refer the reader to [2], [10], [11] and [15]. Using a weak corollary of Theorem 1 of this article (also,
Theorem 1 of [9]), we have explicitly constructed trial wave functions for the {\it minimal IQL configurations} of $N$ Fermions in 
a Jain state with filling factor $ < 1 /2$ (see [11]); it is not possible to give a full account of our recent results here. 
The central result of this article (Theorem 1), presents a useful sufficient condition on a multigraph $\Gamma$ that ensures 
nontriviality of $g(\Gamma)$. There is nothing akin to Theorem 1 in the existing literature. Whenever Theorem 1 is applicable 
to even a single member of $C(N, L, d)$, it readily yields an upper bound on $p(N, L, d)$. Our proof of Theorem 1 is purely algebraic 
in nature; so, the edge-function (or the edge-matrix) of a multigraph is of key importance in the proof. In Theorem 1 we consider only 
those multigraphs $\Gamma$ that can be partitioned into two or more sub-multigraphs $\Gamma _{1}, \ldots , \Gamma _{m}$ such that each 
$g(\Gamma _{i})$ is nonzero (in particular, if $\Gamma _{i}$ has no edges) and the {\it inter-edges} between pairs $\Gamma _{i}$, 
$\Gamma _{j}$ are more `dominating' (in a specific way) than the {\it intra-edges} within each $\Gamma _{i}$. Using Theorem 1, we are 
able to construct several infinite families of invariants (including skew-invariants, see Theorem 2) as well as families of $k$-linearly 
independent semi-invariants of a binary $N$-ic form over $k$ (see Theorem 3). At its core, our approach has its source in [1]; this is
very philosophical and hence almost impossible to articulate. In closing, we share our optimism that there is a generalization of 
Theorem 1 yet to be discovered, that will allow construction of all semi-invariants as symmetrized-graph-monomials. 
 
In what follows, $N$ is tacitly assumed to be an integer $\geq 2$, $k$ denotes a field and 
$z_{1}, \dots, z_{N}$ are indeterminates. We let $z$ stand either for $(z_{1}, \dots , z_{N})$ or the set 
$\{z_{1}, \ldots, z_{N} \}$. It is tacitly assumed that either $k$ has characteristic $0$ or 
the characteristic of $k$ is $ > N$. As usual, given a positive integer $n$, $S_{n}$ denotes the group of all 
permutations of the set $\{1 , \ldots, n \}$. \\

\noindent \underline {\bf Definitions}: Let $m$ and $n$ be positive integers.
\begin{enumerate}
\item  Let $Symm_{N}:\, k [z] \rightarrow k [z]$ be the 
{\it Symmetrization operator} defined by
\[ Symm_{N} (f) := \sum_{\sigma \in S_{N}} f(z_{\sigma (1)}, \ldots , z_{\sigma (N)}) .\]
$f \in k [z]$ is said to be {\it symmetric} provided
\[ f(z_{\sigma (1)}, \ldots , z_{\sigma (N)}) \, = \, f(z_{1}, \ldots , z_{N}) \;\;\; 
\mbox{ for all $\sigma \in S_{N}$. } \]
\item For an $m \times n$ matrix $A := [a_{ij}]$, let $r_{i} (A) := a_{i1}+ \cdots + a_{in}$ 
(the sum of the entries in the $i$-th row of $A$) for $1 \leq i \leq m$ and let
\[ \| A \| \; :=  \; r_{1}(A) + \cdots + r_{m}(A) \; = \; \sum _{i=1}^{m} \sum _{j=1}^{n} a_{ij}. \]
\item Let $E (N)$ denote the set of all $N \times N$ symmetric matrices $A:=[a_{ij}]$ such that
each $a_{ij}$ is a nonnegative integer and $a_{ii} = 0$ for $1 \leq i \leq N$.
\item Given an integer $d$, by $E(N, \, d)$ we denote the subset of $A \in E (N)$ such that 
$r_{i} (A) = d$ for $1 \leq i \leq N$, {\it i.e.}, each row-sum of $A$ is exactly $d$.
\item For an $N \times N$ matrix $A := [a_{ij}]$, let
\[ \delta (z,\, A) \, := \, \prod _{1 \leq i < j \leq N} \, (z_{i} - z_{j})^{a_{ij}} . \]
\item Let $D_{(m, n)} := [(c_{ij}]$ be the $m \times n$ matrix such that
\[ c_{ii} \, :=  \, \left \{ \begin{array}{ll} 0 & \mbox{if $i = j$,} \vspace{0.1in} \\
1 & \mbox{if $i \neq j$.} \end{array} \right . \]
By $D_{n}$, we mean $D_{(n, n)}$. In particular, $D_{1} = 0$. \\
\end{enumerate}

\noindent \underline {\bf Lemma 1}: Let $n$ be a positive integer. For $1 \leq i \leq n$, let 
$g _{i} \in {\Bbb Q}(z)$. Then $g_{1}^2 + g_{2}^2 + \cdots + g_{n}^2 = 0$ if and only if 
$g_{i} = 0$ for $1 \leq i \leq n$. In particular, given a $0 \neq g \in {\Bbb Q}(z_{1}, \dots , z_{N})$ 
and a nonempty subset $S \subseteq S_{N}$, we have
\[ \sum _{\sigma \in S} g (z_{\sigma (1)}, \dots , z_{\sigma (N)})^{2} \, \neq  \, 0 . \] 

\noindent \underline{\bf Proof}: With the notation of (i), assume that $g_{1} \neq 0$. 
Let $h := g_{1}^2 + g_{2}^2 + \cdots + g_{n}^2$. For $1 \leq  i \leq n$, let 
$p_{i} , q_{i} \in {\Bbb Q}[z_{1}, \dots , z_{N}]$ be polynomials such that $g_{i} q_{i} = p_{i}$ 
and $q_{i} \neq 0$. Note that, $g_{1} \neq 0$ implies $p_{1} \neq 0$. Now since 
$f:= p_{1} q_{1} q_{2} \cdots q_{n}$ is a nonzero polynomial with coefficients in ${\Bbb Q}$, there 
exists $(a_{1}, \dots , a_{N}) \in {\Bbb Q}^{N}$ such that $f(a_{1}, \dots , a_{N}) \neq 0$. Fix such 
an $N$-tuple $(a_{1}, \dots , a_{N})$ and let $c_{i} := g_{i} (a_{1}, \dots , a_{N})$ for $1 \leq i \leq n$. 
Then, $c_{1} \neq 0$ and $c_{i} \in {\Bbb Q}$ for $1 \leq i \leq n$. Since $c_{1}^2 > 0$ and 
$(c_{2}^2 + \cdots + c_{n}^2) \geq 0$, we have $h (a_{1}, \dots , a_{N}) > 0$. This proves the first 
claim of (i); the second claim of (i) easily follows. Assertion (ii) readily follows from (i). $\Box$ \\

\noindent \underline {\bf Definitions}: 
\begin{enumerate}
\item For $B \subseteq \{1, 2, \dots , N \}$, let
\[ \pi (B) \, := \, \{ (i, j) \in B \times B \, \mid \, i < j \}.  \]
By abuse of notation, $\pi (B)$ is also identified as the set of all $2$-element
subsets of $B$. The set $ \pi (\{1, \dots , N \})$ is denoted by $\pi [N]$.
\item Given $C \subseteq \pi [N]$ and a function $\varepsilon : C \rightarrow {\Bbb N}$, the image of 
$(i, j) \in C$ via $\varepsilon $ is denoted by $\varepsilon (i, j)$. An integer $w \in {\Bbb N}$ is 
identified with the constant function $C \rightarrow {\Bbb N}$ such that $(i, j) \rightarrow w$ for all
$(i, j) \in C$.
\item Given $C \subseteq \pi [N]$ and a function $\varepsilon : C \rightarrow {\Bbb N}$, define
\[ v(z, C, \varepsilon) \, := \, \prod_{(i, j) \in C} \, (z_{i} - z_{j})^{\varepsilon (i, j)} \]
with the understanding that $v(z, \emptyset, \varepsilon) = 1$.  \\
\end{enumerate}

\noindent \underline {\bf Remark 1}:
There is an obvious bijective correspondence $\varepsilon \leftrightarrow [a_{ij}]$ between the set 
of functions $\varepsilon : \pi [N] \rightarrow {\Bbb N}$ and the set $E (N)$, given by
\[ a_{ij} \,= \, \varepsilon (i, j)\;\;\;\mbox{for $1 \leq i < j \leq N$.} \] \\

Suppose $m_{1} \leq m_{2} \leq \cdots \leq m_{q}$ is a partition of $N$ and $M \in E(N)$. Consider $M$ as a 
$q \times q$ block-matrix $[M_{rs}]$, where $M_{rs}$ has size $m_{r} \times m_{s}$ for $1 \leq r, s \leq q$. 
View $M$ as the sum $M^{*} + M^{**}$, where $M^{*}$ is the $q \times q$ block-diagonal matrix having $M_{rr}$ 
as its $r$-th diagonal block and where $M^{**}$ is the $q \times q$ block-matrix whose diagonal blocks are 
zero-matrices. Clearly, $M^{*}$ and $M^{**}$ both are in $E(N)$ and $M_{rr} \in E(m_{r})$ for $1 \leq r \leq q$. \\

\noindent \underline {\bf Definitions}: Let the notation be as above.
\begin{enumerate}
\item For $1 \leq r \leq q$, define
\[ A_{r} \; := \; \left \{ i + m_{0} + \cdots + m_{r - 1}   \, \mid  \, 1 \leq i \leq m_{r}  \right \} . \]
\item For $1 \leq r \leq q$, let $G_{r}$ denote the group of permutations of the set $A_{r}$.
\item Define
\[ \pi \; := \; \bigcup _{1 \leq r < s \leq q} A_{r} \times A_{s}. \]
\item For $1 \leq r \leq q$ and $(i, j) \in \pi (A_{r})$, let $\varepsilon _{r} (i, j)$ denote the $ij$-th
entry of $M^{*}$.
\item For $1 \leq r \leq q$, define
\[ \delta _{r} ( M^{*} ) \; := \; Symm _{m_{r}} \left( v(z, \pi (A_{r}), \varepsilon_{r}) \right ). \]
\item For $(i, j) \in \pi [N]$, let $\varepsilon (i, j)$ denote the $ij$-th entry of $M^{**}$. \\
\end{enumerate}

\noindent \underline {\bf Remark 2}:
\begin{enumerate}
\item Observe that
\[ \pi \; = \; \pi [N] \setminus \bigcup _{i = 1}^{q} \pi (A_{i}). \]
\item For each $r$, the $\varepsilon _{r} (i, j)$ are the the entries in the strict upper-triangle of the 
symmetric matrix $M_{rr}$. 
\item We have $\delta (z, M^{**}) \; = \; v(z, \pi [N], \varepsilon)$ and
\[\delta (Z, M^{*}) \; = \;  \prod _{r = 1}^{q} v(z, \pi (A_{r}), \varepsilon_{r}) .\]
\item We have $\delta (z, M) = \delta (z, M^{*}) \cdot \delta (z, M^{**})$.
\item For each $r$, we have
\[\delta _{r} ( M^{*} ) \; = \; \sum _{\sigma \in G_{r}} \sigma (v(z, \pi (A_{r}), \varepsilon_{r})). \]
\item The $\varepsilon (i, j)$ are the entries in the strict upper-triangle of the symmetric matrix $M^{**}$. \\
\end{enumerate}

\noindent \underline {\bf Theorem 1}: Let the notation be as above. Assume $q \geq 2$ and of the following 
properties (1) - (3), either  (1) and (2) hold or  (1) and (3) hold.
\begin{description}
\item[(1)] For $1 \leq r < s \leq q$, the matrix $M_{rs}$ has only  positive entries.
\item[(2)] For  $1 \leq r < s \leq q$, the positive integer  $b (m_{r}, m_{s}) \, := \, \| M_{rs} \|$ depends
only on the ordered pair $(m_{r}, \, m_{s})$ and furthermore, if $m_{r} = m_{s}$, then $b (m_{r}, m_{s})$ 
is an even integer.
\item[(3)] Characteristic of $k$ is $0$ and for $1 \leq r < s \leq q$, $\| M_{rs} \|$ is even.
\end{description}
Also, assume that the properties (i) - (iv) listed below are satisfied.
\begin{description}
\item[(i)] Either $m_{i} < m_{j}$ for $1 \leq i < j \leq q$ or $M^{*} = 0$.
\item[(ii)] If properties (1) and (2) hold, then $\prod _{r = 1}^{q} \delta _{r} ( M^{*} ) \neq 0$.
\item[(iii)] If property (2) does not hold but properties (1) and (3) hold, then each entry of 
$M^{*}$ is an even integer.
\item[(iv)] The least nonzero entry of the matrix $M^{**}$ is strictly greater than the greatest entry 
of the matrix $M^{*}$.
\end{description}
Then $Symm_{N} \left (  \delta (z, M) \right ) \; \neq \; 0$. \\

\noindent \underline{\bf Proof}: Define $m_{0} = 0$. At the outset, observe that a permutation
$\sigma \in S_{N}$ can be naturally viewed as a permutation of $\pi [N]$
by letting $\sigma (i, j) \, := \, \{ \sigma (i), \sigma (j) \}$, {\it i.e.},
for $(i, j) \in \pi [N]$,
\[\sigma (i, j) \, := \,
\left \{\begin{array}{ll} (\sigma (i), \sigma (j)) &
\mbox{if $\sigma (i) < \sigma (j)$,} \vspace{0.1in} \\ (\sigma (j), \sigma (i)) &
\mbox{if $\sigma (j) < \sigma (i)$.} \end{array} \right . \]
Thus $S_{N}$ is regarded as a subgroup of the group of permutations of $\pi [N]$. 

For $\sigma \in S_{N}$ and $1 \leq r \leq q$, define
\[ B_{r} (\sigma)\, := \, \sigma ^{-1} (A_{r}) \, = \,
\{ i \, \mid \, 1 \leq i \leq N \;\;\mbox{and}\;\;\sigma (i) \in A_{r} \}. \]
Clearly, sets $B_{1} (\sigma), \dots , B_{q} (\sigma)$ partition $\{1, \dots , N \}$ and $B_{i}$
has cardinality $m_{i}$ for all $1 \leq i \leq q$.

Define 
\[ G \, := \, \{ \sigma \in S_{N} \, \mid \, \sigma (i, j) \in \pi
\;\;\mbox{for all $(i, j) \in \pi$} \} . \]
For $1 \leq r \leq q$, a permutation $\sigma \in G_{r}$ is to be regarded as an element of $S_{N}$ by
declaring $\sigma (i) = i$ if $i \in \{1, \ldots , N \}\setminus A_{r}$. This way each $G_{r}$ is 
identified as a subgroup of $S_{N}$.  

Given $\sigma \in G$ and $(i, j) \in \pi (A_{r})$ with $1 \leq r \leq q$,
clearly there is a unique $s$ with $1 \leq s \leq q$ such that $\sigma (i , j) \in \pi (A_{s})$.
Fix a $\sigma \in G$. Consider $i \in B_{r} (\sigma) \cap A_{s}$ with $1 \leq s \leq q$.
Then for $i \neq j \in A_{s}$, we must have $\{\sigma (i), \sigma (j) \}$ in $\pi (A_{r})$
and hence $j \in B_{r} (\sigma)$. It follows that $A_{s} \subseteq B_{r} (\sigma)$. If
$1 \leq s < p \leq q$ are such that $A_{s} \cup A_{p} \subseteq B_{r} (\sigma)$, then
an $(i, j) \in A_{s} \times A_{p}$ is in $\pi$ whereas $\sigma (i, j)$ is
in $\pi (A_{r})$. This is impossible since $\sigma \in G$. Thus we have established
the following: given $r$ with $1 \leq r \leq q$ and $\sigma \in G$, there is a unique
integer $r (\sigma)$ such that $1 \leq r(\sigma) \leq q$ and $B_{r} (\sigma) = A_{r (\sigma)}$.
In other words, the image sets $\sigma (A_{1}), \dots , \sigma (A_{q})$ form a
permutation of the sets $A_{1}, \dots , A_{q}$. If $1 \leq r < s \leq q$ and 
$\sigma \in G$, then since $r (\sigma) \neq s (\sigma)$, we infer that
\[ \pi \cap \left ( A_{r (\sigma)} \times A_{s (\sigma)} \right ) \, \neq \, \emptyset
\;\;\;\mbox{if and only if $r (\sigma) < s (\sigma)$.} \]
Moreover,
\[ m_{r (\sigma)} \, = \, m_{r} \;\;\; \mbox{for all $1 \leq r \leq q$ and $\sigma \in G$.} \]
If the first case of (i) holds, {\it i.e.}, the integers $m_{i}$ are mutually unequal, then we must
have $r(\sigma) = r$ for all $1 \leq r \leq q$ and $\sigma \in G$. Hence, in this case $G$ is the 
direct product of (the mutually commuting) subgroups $G_{1}, G_{2}, \dots , G_{q}$.

Hypothesis (1) implies $v(z, \pi[N], \varepsilon) = v(z, \pi, \varepsilon)$. If 
$G  = G_{1} \times G_{2} \times \cdots  \times G_{q}$, then we have
\[ \sum _{\sigma \in G} \left ( \prod _{r = 1}^{q} \sigma (v(z, \pi (A_{r}), \varepsilon _{r})) \right ) \, = \,
\prod _{r=1}^{q} \left ( \sum _{\theta \in G_{r}} \theta (v(z, \pi (A_{r}), \varepsilon _{r})) \right ) .\]
For $1 \leq r \leq q$, define
\[ w_{r} \, := \, \sum _{(i, j) \in \pi (A_{r})} \varepsilon _{r} (i, j) \;\;\;\; \mbox{and} \;\;\;\;
w \, := \, \sum_{i = 1}^{q} w_{i} . \]
Our hypothesis (i) assures that if $m_{i} = m_{j}$ for some $i \neq j$, then $w = 0$.

Now let $t, t_{1}, \ldots , t_{q},  x_{1}, \dots , x_{N}$ be indeterminates and let
\[ \alpha : k [z_{1}, \dots , z_{N}] \rightarrow
k [t, t_{1}, \cdots , t_{q}, x_{1}, \dots , x_{N}] \]
be the injective $k$-homomorphism of rings defined by
\[ \alpha (z_{i}) \, := \, t x_{i} + t_{r} \;\;\;
\mbox{if $i \in A_{r}$ with $1 \leq r \leq q$.}  \]
Then given $\sigma \in S_{N}$, $(i, j) \in \pi [N]$ and $1 \leq r, s \leq q$, we have
\[ \alpha (z_{\sigma (i)} - z_{\sigma (j)}) \, = \,
t( x_{\sigma (i)} - x_{\sigma(j)} ) + (t_{r} - t_{s}) \]
if and only if $(\sigma (i), \sigma (j)) \in A_{r} \times A_{s}$.

Let $x$ stand for $(x_{1}, \dots , x_{N})$ and $T$ stand for $(t_{1}, \dots , t_{q})$.
Given $f \in k [t, T, X]$, by the {\it $x$-degree} (resp. {\it $T$-degree}) of $f$,
we mean the total degree of $f$ in the indeterminates $x_{1}, \dots , x_{N}$
(resp. $t_{1}, \dots , t_{q}$). Now fix a $\sigma \in G$ and consider
\[ V_{\sigma} (x, t, T)\, := \, \alpha ( \sigma (v (z, \pi, \varepsilon ))). \]
For an ordered pair $(i, j)$ with $1 \leq i, j \leq q$, set
\[ A (\sigma, i, j) \, := \, \pi \cap ( A_{i (\sigma)} \times A_{j (\sigma)}) .\]
It is straightforward to verify that $V_{\sigma} (x, 0, T)$ is
\[ \prod _{1 \leq r < s \leq q} \left ( \prod _{(i, j) \in A (\sigma, r, s)} (t_{r} - t_{s})^{\varepsilon (i, j)}
  \cdot \prod _{(i, j) \in A (\sigma, s, r)} (t_{s} - t_{r})^{\varepsilon (i, j)} \right ).   \]
Suppose condition (2) of the theorem holds. Then for $1 \leq r < s \leq q$, we have
\[\sum _{(i, j) \in A (\sigma, r, s)} \varepsilon (i, j) \, = \,
\left \{\begin{array}{ll} 0 & \mbox{if $s (\sigma) < r (\sigma)$,} \vspace{0.1in} \\
b (m_{r}, m_{s}) & \mbox{if $r (\sigma) < s (\sigma)$.} \end{array} \right . \]
Further, if $1 \leq r < s \leq q$ are such that $s (\sigma) < r (\sigma)$, then
\[ m_{s} = m_{s (\sigma)} \leq m_{r(\sigma)} = m_{r} \;\;\;\mbox{implies 
$m_{s} = m_{s (\sigma)}  = m_{r(\sigma)} = m_{r}$} \]
and so, (2) assures that $b (m_{r}, m_{s})$ is an even integer. Hence, if property (2) holds, then  
\[ V_{\sigma} (x, 0, T)\, := \, \prod _{1 \leq r < s \leq q} (t_{r} - t_{s})^{ b(m_{r}, m_{s})} .\]
On the other hand, if condition (3) holds, then we merely observe that there is a nonzero homogeneous 
$g_{\sigma} \in {\Bbb Q}[t_{1}, \dots , t_{q}]$ such that $V_{\sigma} (x, 0, T) = g_{\sigma}^{2}$.
In any case, the $t$-order of $V_{\sigma} (x, 0, T)$ is $0$ ({\it i.e.}, $V_{\sigma} (x, t, T)$ is not
a multiple of $t$) and the $T$-degree of $V_{\sigma} (x, 0, T)$ is
\[ d \, := \, \sum _{(i, j) \in \pi} \varepsilon (i, j). \]
Define
\[ \gamma\, := \, \sum _{\sigma \in G} \, \sigma ( v(z, \pi, \varepsilon) ) \;\;\;\mbox{and}
\;\;\; V (x, t, T) \, := \, \sum _{\sigma \in G}\, V_{\sigma} (x, t, T) . \]
Then $\alpha (\gamma) =  V (x, t, T)$. If (2) holds, then letting $|G|$ denote the cardinality of $G$,
we have $|G| \neq 0$ in $k$ and
\[ (\#)\;\;\;\;\;\; V (x, 0, T) \, = \, |G|\,\prod _{1 \leq r < s \leq q} (t_{r} - t_{s})^{b (m_{r}, m_{s})} \]
and hence $V (x, 0, T) \neq 0$. On the other hand, if (3) holds, then we have
\[ V(x, 0, T) \, = \, \sum _{\sigma \in G}\, g_{\sigma}^{2}, \]
which is necessarily nonzero in view of Lemma 1. Now it is clear that
$\alpha (\gamma) \neq 0$, the $t$-order of $\alpha (\gamma)$ is $0$ and the $T$-degree of
$\alpha (\gamma)$ is $d$.

For $\sigma \in S_{N}$, define
\[F_{\sigma} (z) \, := \, \prod _{r = 1}^{q} \sigma (v(z, \pi (A_{r}), \varepsilon _{r})) \;\;\;\;
\mbox{and} \;\;\;\;
 W_{\sigma} (x, t, T) \, := \, 
 \prod _{r = 1}^{q} \alpha (\sigma (v(z, \pi (A_{r}), \varepsilon _{r}))). \]
Then $W_{\sigma} (x, t, T)= \alpha (F_{\sigma} (z))$, If $\varepsilon_{r} = 0$ for all $r$, then $F_{\sigma} (z) = 1$ 
and hence 
\[ \sum_{\sigma \in G} F_{\sigma} (x) \, = \, |G| \, \neq \, 0. \]
If $G = G_{1} \times \cdots \times G_{q}$, then we have
\[\sum_{\sigma \in G} F_{\sigma} (x) \, = \, 
\prod _{r = 1}^{q} \left ( \sum _{\theta \in G_{r}} \theta (v(z, \pi (A_{r}), \varepsilon _{r})) \right ) .\]
Now suppose $G = G_{1} \times \cdots \times G_{q}$.
Given $\sigma \in G$, write $\sigma =: \theta _{1} \theta _{2} \cdots \theta _{q}$, where $\theta _{r} \in G_{r}$ 
for $1 \leq r \leq q$. Then
\[ \alpha (\sigma (v(z, \pi (A_{r}), \varepsilon _{r})))\, = \,
t^{w_{r}}\, \theta _{r} (v(x, \pi (A_{r}), \varepsilon _{r}))\, = \,t^{w_{r}}\, \sigma (v(x, \pi (A_{r}), \varepsilon _{r}))  \]
and hence
\[ W_{\sigma} (x, t, T) \, = \, t^{w} \, \prod _{r = 1}^{q} \sigma (v(x, \pi (A_{r}), \varepsilon _{r}))
\, = \, t^{w} F_{\sigma} (x). \]
Consequently, 
\[\alpha (\sigma (v(z, \pi, \varepsilon)))
\prod _{r = 1}^{q} \alpha (\sigma ( v(z, \pi (A_{r}), \varepsilon _{r}) )) \, = \,
t^{w} V_{\sigma}(x, t, T) F_{\sigma} (x). \]

Case I: hypothesis (ii) holds. Then as proved above $V_{\sigma}(x, 0, T)$ is independent of the choice of 
$\sigma \in G$ and $V_{\sigma}(x, 0, T)$ is a nonzero polynomial depending only on $T$. In particular, letting 
$\iota \in S_{N}$ denote the identity permutation, we have $V_{\iota} (x, 0, T) \neq 0$ and
\[ \sum_{\sigma \in G} V_{\sigma} (x, 0, T) F_{\sigma} (x) \, = \,
V_{\iota} (x, 0, T) \sum _{\sigma \in G} F_{\sigma} (x). \]
The sum appearing on the right of the above equation is obviously independent of $t$; moreover, hypothesis (ii) 
assures that it is nonzero and thus has $t$-order $0$. Case II: hypothesis (iii) holds. Then 
$V_{\sigma} (x, 0, T) = g_{\sigma}^{2}$ as well as $F_{\sigma} (x) = f_{\sigma}^{2}$, where $g_{\sigma} \in k[T]$ 
and $f_{\sigma} \in k[x]$ are nonzero polynomials. In this case, Lemma 1 assures that
\[\sum_{\sigma \in G} V_{\sigma} (x, 0, T) F_{\sigma} (x) \, = \,
\sum_{\sigma \in G} (f_{\sigma} g_{\sigma})^{2} \, \neq \, 0. \]
In either case, the sum
\[\sum _{\sigma \in G} V_{\sigma}(x, t, T) W_{\sigma} (x, t, T) \, = \, 
\sum _{\sigma \in G} t^{w} V_{\sigma}(x, t, T) F_{\sigma} (x) \]
has $t$-order exactly $w$. 

Next, for $\sigma \in S_{N}$, let
\[ R (\sigma) \, := \, \bigcup _{1 \leq r \leq q} \pi (B_{r} (\sigma)) . \]
Observe that $ \pi \cap R (\sigma) =  \emptyset $ if and only if $\sigma \in G$. Also, observe that
\[ \alpha (z_{\sigma (i)} - z_{\sigma (j)}) \, = \, t ( x_{\sigma (i)} - x_{\sigma (j)} ) + (t_{r} - t_{s}), \]
where $r = s$ if and only if $(i, j) \in R (\sigma)$. 

Fix a $\sigma \in S_{N} \setminus G$. Then clearly
\[v(z, \pi, \varepsilon) \, = \, v(z, \pi[N], \varepsilon) \, = \, 
v(z, R (\sigma), \varepsilon) v (z, \pi [N] \setminus R (\sigma), \varepsilon). \]
Moreover, note that
\[v(z, R (\sigma), \varepsilon) \, = \, v (z, \pi \cap R (\sigma), \varepsilon) \;\;\;\; \mbox{and} \;\;\;\; 
v(z, \pi [N] \setminus R (\sigma), \varepsilon) \, = \, v (z, \pi \setminus R (\sigma), \varepsilon). \]
Define
\[ \lambda (\sigma) \, := \, \sum _{(i, j) \in \pi \cap R (\sigma)} \varepsilon (i, j)
\;\;\;\; \mbox{and} \;\;\;\; 
d (\sigma) \, := \, \sum _{(i, j) \in \pi \setminus R (\sigma)} \varepsilon (i, j) . \]
Then $d (\sigma) \, = \, d - \lambda (\sigma)$. From our choice of $\sigma$ and hypothesis (1), it follows
that $\lambda (\sigma) \geq 1$ and hence $d (\sigma) < d $. Let
\[P_{\sigma} (x, t, T) \, := \, \alpha (\sigma (v (z, \pi \cap R (\sigma), \varepsilon))), \;\;\;\;
 Q_{\sigma} (x, t, T) \,:= \,\alpha (\sigma (v (z, \pi \setminus R (\sigma), \varepsilon))).  \]
Observe that $V_{\sigma}(x, t, T) = P_{\sigma} (x, t, T) \cdot Q_{\sigma} (x, t, T)$,
\[ P_{\sigma} (x, t, T) \, = \, t^{\lambda (\sigma)} \cdot
\prod _{(i, j) \in \pi \cap R (\sigma)} (x_{\sigma (i)} - x_{\sigma (j)})^{\varepsilon (i, j)} \]
and $Q_{\sigma}(x, 0, T)$ is a nonzero $T$-homogeneous polynomial of $T$-degree $d (\sigma)$. Hence
the $t$-order of $V_{\sigma}(x, t, T)$ is exactly $\lambda (\sigma)$. For $1 \leq r \leq q$, let
\[\begin{array}{lll} P_{\sigma}^{(r)}(x, t, T) & := &  \alpha (\sigma (v (z, \pi (A_{r}) \cap R (\sigma), \varepsilon _{r}))) , \\
Q_{\sigma}^{(r)} (x, t, T)  & :=  & \alpha (\sigma (v (z, \pi (A_{r}) \setminus R (\sigma), \varepsilon_{r}))). \end{array}  \]
Now for $1 \leq r \leq q$, we do have
\[\sigma (v(z, \pi(A_{r}), \varepsilon _{r})) \, = \, \sigma (v (z, \pi (A_{r}) \cap R (\sigma), \varepsilon _{r})) \cdot
\sigma (v (z, \pi (A_{r}) \setminus R (\sigma), \varepsilon_{r})) \]
and hence 
\[\alpha (\sigma (v(z, \pi(A_{r}), \varepsilon _{r}))) \, = \, P_{\sigma}^{(r)}(x, t, T) \cdot Q_{\sigma}^{(r)} (x, t, T) .\]
Since $\pi (B_{s} (\sigma)) \cap \pi (B_{r} (\sigma) = \emptyset = \pi (A_{r}) \cap \pi (A_{s})$ for $1 \leq r < s \leq q$,
we have
\[ \pi \cap R (\sigma) \, = \, \{(i, j) \in \pi \, \mid \, \sigma (i, j) \in \pi [N] \setminus \pi \} \, = \,
\bigsqcup _{r = 1}^{q}  \left ( \pi \cap \pi (B_{r} (\sigma)) \right )\]
and
\[ J \, := \, \bigsqcup _{r = 1}^{q} \left ( \pi (A_{r}) \setminus R (\sigma) \right ) \, = \, 
\{(i, j) \in \pi [N] \setminus \pi \, \mid \, \sigma (i, j) \in \pi \}. \]
Recall that $\sigma$ is also viewed as a permutation of $\pi [N]$. Hence $J$ and $\pi \cap R(\sigma)$ have the
same cardinality. Partition $\pi \cap R (\sigma)$ into $q$ subsets $I_{1} (\sigma), \dots , I_{q} (\sigma)$ such that 
$| I_{r} (\sigma) | = | \pi (A_{r}) \setminus R (\sigma) |$ for $1 \leq r \leq q$. For $1 \leq r \leq q$, define 
\[\lambda _{r} (\sigma) \, := \, \sum _{(i, j) \in I_{r} (\sigma)} \varepsilon (i, j) \;\;\;\;
\mbox{and} \;\;\;\; e_{r} (\sigma) \, := \, \sum _{(i, j) \in \pi (A_{r}) \cap R (\sigma)} \varepsilon_{r} (i, j). \]
Then $\lambda (\sigma) = \lambda _{1} (\sigma) + \cdots + \lambda _{q} (\sigma)$, the $t$-order of 
$P_{\sigma}^{(r)}(x, t, T)$ is $e_{r} (\sigma)$ and the $t$-order of $Q_{\sigma}^{(r)} (x, t, T)$ is $0$ for 
$1 \leq r \leq q$. Consequently, the $t$-order of $V_{\sigma} (x, t, T) W_{\sigma} (x, t, T)$ is 
\[ \lambda (\sigma) + \sum _{r = 1}^{q} e_{r} (\sigma) \, = \, 
\sum _{r = 1}^{q}  e_{r} (\sigma) + \lambda _{r} (\sigma) . \]
Our hypothesis (iv) guarantees that firstly $e_{r} (\sigma) + \lambda _{r} (\sigma) \geq w_{r}$ for $1 \leq r \leq q$ 
and secondly, since $\sigma$ is not in $G$, there is at least one $r$ with 
$e_{r} (\sigma) + \lambda _{r} (\sigma) \geq  w_{r} + 1$. It follows that for each
$\sigma \in S_{N} \setminus G$, the $t$-order of $V_{\sigma} (x, t, T) W_{\sigma} (x, t, T)$ is at least $w + 1$.

Let $\Upsilon := Symm_{N} (\delta(z, M))$. Then we have
\[ \Upsilon \,  = \, Symm_{N} \left ( v(z, \pi, \varepsilon)
\prod _{r = 1}^{q} v(z, \pi (A_{r}), \varepsilon _{r}) \right) \]
and hence
\[ \alpha (\Upsilon) \, = \,\sum _{\sigma \in G} V_{\sigma}(x, t, T) W_{\sigma} (x, t, T) +
\sum _{\sigma \in G \setminus S_{N}} V_{\sigma }(x, t, T) W_{\sigma } (x, t, T) .\]
Since $G$ is nonempty, the first sum on the right of the above equality is nonzero. From what has been shown above
the first sum on the right has $t$-order $w$ whereas the second sum on the right has $t$-order at least $w + 1$. Hence
$\alpha (\Upsilon)$ has $t$-order $w$. Since $w$ is a nonnegative integer, $\alpha (\Upsilon) \neq 0$. In particular, 
$\Upsilon \neq 0$. $\Box$ \\

\noindent \underline {\bf Remark 3}: We continue to use the above notation.
\begin{enumerate}
\item Suppose $M$ satisfies the hypotheses of Theorem 1 and $\lambda$ is a positive integer such that 
$Symm _{m_{r}} ( \delta (z, \lambda M_{rr})) \neq 0$ for $1 \leq r \leq q$. Then $\lambda M$ also satisfies
the hypotheses of Theorem 1. In general, the polynomials $Symm_{N} (\delta (z, M))$ and $Symm_{N} (\delta (z, \lambda M))$ 
do not seem to be related in any obvious manner (see the last of the Examples 1 below).
\item Suppose for $1 \leq i \leq s$, there is a partition $m^{(i)}$ of $N$ with respect to which $M _{i} \in E(N)$ satisfies 
the hypotheses of Theorem 1 and let $\Upsilon _{i}:= Symm_{N} (\delta(z, M_{i}))$. If 
$\alpha (\Upsilon _{1}), \ldots , \alpha (\Upsilon _{s})$ are $k$-linearly independent, then 
$\Upsilon _{1}, \ldots , \Upsilon _{s}$ are also $k$-linearly independent. Now to ensure 
$k$-linear independence of $\alpha (\Upsilon _{1}), \ldots , \alpha (\Upsilon _{s})$, it suffices to ensure the $k$-linear 
independence of their respective $t$-initial forms. For simplicity, assume that property (2) is satisfied by the $M_{i}$ and 
$M_{i}^{*} = 0$ for $1 \leq i \leq s$. Then from the equality $(\#)$ in the proof of Theorem 1 it follows that the $t$-initial 
coefficient,{\it i.e.}, the coefficient of the lowest power of $t$ present, of each $\alpha (\Upsilon _{i})$ is of the type 
$ c \prod _{1 \leq r < s \leq q} (t_{r} - t_{s})^{b(m_{r}, m_{s})}$ for some $0 \neq c \in k$. The $k$-linear independence of 
such products is completely determined by the exponents $b(m_{r}, m_{s})$. \\
\end{enumerate}

\noindent \underline {\bf Example 1}: 
\begin{enumerate}
\item Consider the following $E_{1}, E_{2}, E_{3} \in E(6)$ presented as $2 \times 2$ block-matrices.
\[E_{i} \, := \, \mat{2}{0 & C_{i} \\ C_{i}^{T} & 0}, \]
where
\[ C_{1} \, := \, \mat{3}{3 & 3 & 3 \\ 3 & 3 & 3 \\ 3 & 3 & 4}, \;\;\; 
C_{2} \, := \, \mat{3}{3 & 3 & 3 \\ 3 & 4 & 3 \\ 3 & 3 & 4}, \;\;\;
C_{3} \, := \, \mat{3}{3 & 3 & 3 \\ 3 & 3 & 4 \\ 3 & 3 & 4} . \]
A direct computation using MAPLE shows that
$Symm_{6} ( \delta (z, E_{1})) \neq 0$, $Symm_{6} ( \delta (z, E_{2})) = 0$ and $Symm_{6} ( \delta (z, E_{3})) \neq 0$. 
Of course, in the case of $E_{1}$, Theorem 1 does apply. Since $\| C_{2} \| = 29 = \| C_{3} \|$ is an odd integer, 
Theorem 1 can not be applied in the case of $E_{2}$, $E_{3}$. 
\item For $j = 1, 2$, let $E_{j} \in E(5, 18)$ be presented in $2 \times 2$ block-format as
\[ E_{j} \, := \, \mat{2}{0 & A_{j} \\ A_{j}^{T} & B}, \;\;\;\mbox{where}\;\;\;
B \, := \, \mat{3}{0 & 1 & 7 \\ 1 & 0 & 1 \\ 7 & 1 & 0}, \]
\[ A_{1} \, := \, \mat{3}{5 & 13 & 0 \\ 5 & 3 & 10} \;\;\;\mbox{and}\;\;\;
A_{2} \, := \, \mat{3}{8 & 10 & 0 \\ 2 & 6 & 10}. \]
Then a MAPLE computation shows that $h_{j} := Symm_{5} (\delta (z, E_{j})) \neq 0 $ for $j = 1, 2$. 
Up to a nonzero integer multiple, $h_{1}$ and $h_{2}$ are the same; either one can be identified as the Hermite's 
invariant of a quintic binary form (see [2] or [3]). Since this invariant has weight $45$, it is a skew invariant. Let 
$M \in E(9, 90)$ be the $2 \times 2$ block-matrix $[M_{ij}]$ such that $M_{11} = 0$, $M_{12}$ is the $4 \times 5$ 
matrix having each entry $18$ and $M_{22} \in \{E_{1}, E_{2} \}$. Note that Theorem 1 is applicable and thus 
$g := Symm_{9} (\delta (z, M))$ is a nonzero invariant of a binary nonic. Also, since $g$ has weight $405$, $g$ is 
a skew invariant. 
\item Let $M \in E(4, 2)$ be the $2 \times 2$ block matrix $[M_{ij}]$, where 
$M_{11} = 2 D_{2} = M_{22}$ and $M_{12} = 0 = M_{21}$. Let $g := Symm _{4} (\delta (z, M))$ and 
$h := Symm _{4} (\delta (z, 2 M))$. Then $2 M \in E(4, 4)$ and by Lemma 1, $g h \neq 0$. Clearly, $g$ and $h$ both are 
invariants of a binary quartic. A computation employing MAPLE shows that $g$ and $h$ are algebraically independent over $k$. \\
\end{enumerate}

\noindent \underline {\bf Lemma 2}: Suppose $d$ is a positive integer such that $N d$ is an integer
multiple of $4$. Then there is an explicitly described $E \in E(N, d)$ such that each entry of $E$ 
is an even integer. Moreover, if $k$ has characteristic $0$, then $g := Symm _{N} (\delta (z, E))$ is 
a nonzero invariant (of degree $d$) of a binary form of degree $N$. \\

\noindent \underline {\bf Proof}: First, suppose $N = 2 m$ for some positive integer $m$ and $d$
is an even positive integer. Let $E \in E(N)$ be the $m \times m$ block matrix $[M_{ij}]$ such that 
$M_{rr} := d D_{2}$ for $1 \leq r \leq m$ and $M_{ij} = 0$ for $1 \leq i < j \leq m$. Then clearly 
$E \in E(N, d)$ and since $d$ is even, each entry of $E$ is an even integer. Secondly, suppose $N$ is odd
and $d = 4 e$ for some positive integer $e$. Our construction proceeds by induction on $N$. If $N = 3$,
then let $E := (2 e) D_{3}$. Henceforth, assume $N \geq 5$. If $N - 3$ is odd, then by induction hypothesis,
we have an $M \in E(N - 3, d)$ such that each entry of $M$ is an even integer. If $N - 3$ is even, then by 
the first part of our proof we have an $M \in E(N - 3, d)$ such that each entry of $M$ is an even integer.
Now let $E$ be the $2 \times 2$ block matrix $[C_{ij}]$ with $C_{11} := (2 e) D_{3}$, $C_{22} := M$ and
$C_{12} = 0 = C_{21}$. Then clearly $E \in E(N, d)$ and each entry of $E$ is an even integer. In either 
case, provided $char\,k = 0$, Lemma 1 assures that $g \neq 0$. $\Box$ \\

\noindent \underline {\bf Theorem 2}: Assume that $N \geq 3$.
\begin{description}
\item[(i)] Suppose $m$, $n$ are positive integers such that $n \geq 2$ and $N = m n$. Let $a$, $b$ 
be positive integers and let $d := 2 a (n - 1) + (m - 1) (n - 1) b$. Then there is an explicitly 
described $E \in E(N, d)$ such that $g := Symm _{N} (\delta (z, E))$ is a (degree $d$) nonzero invariant 
of a binary form of degree $N$.
\item[(ii)] Suppose $m$, $n$, $r$ are positive integers such that $n \geq 2$, $1 \leq r \leq m n - 1$ and 
$N = 2 m n - r$. Given positive integers $a$, $b$ such that 
\[ c \; := \; \frac{2 (n - 1) a + (m - 1) (n  - 1) b}{r} \;\;\; \;\mbox{is an integer,} \]
there is an explicitly described $E \in E(N, m n c)$ yielding a (degree $m n c$) nonzero invariant 
$g := Symm _{N} (\delta (z, E))$ of a binary form of degree $N$. 
\item[(iii)] Suppose $l$, $m$, $n$ are positive integers such that $l < m < n < l + m$ and $N = l + m + n$. Given
a positive integer $d$ such that each of
\[a \, := \, \frac{(m + l - n) d}{ 2 l m}, \;\;\; b \, := \, \frac{(l + n - m) d}{ 2 l n}, \;\;\;
c \, := \, \frac{(m + n - l) d}{ 2 m n} \]
is an integer, there is an explicitly described $E \in E(N, d)$ yielding a (degree $d$) nonzero invariant 
$g := Symm _{N} (\delta (z, E))$ of a binary form of degree $N$. 
\item[(iv)] Suppose $s$ is a nonnegative integer and $t$, $u$, $v$ are positive integers such that 
$t \leq 2 u \leq 2 t - 1$. Then letting 
\[N \, := \, 2 (2 t v + 1) \;\;\;\;\mbox{and} \;\;\;\; d \,:= \,(2 s + 1) (2 u + 1) (4 u v + 2 v + 1), \] 
there is an explicitly described $E \in E(N, d)$ such that $g := Symm _{N} (\delta (z, E))$ is a nonzero invariant 
of a binary form of degree $N$. Moreover, $g$ is a skew invariant of weight 
$w := (2 s + 1) (2 t v + 1) (2 u + 1) (4 u v + 2 v + 1)$. 
\item[(v)] Given $E \in E(N, d)$ such that each entry of $E$ is strictly less than $d$ and 
$Symm _{N} (\delta (z, E)) \neq 0$, a matrix $E^{*} \in E(2 N - 1, d N)$ can be so constructed that 
$g := Symm _{N} (\delta (z, E^{*}))$ is a nonzero invariant of a binary form of degree $2 N - 1$. \\
\end{description}

\noindent \underline {\bf Proof}: To prove (i), let $E \in E(N)$ be the $n \times n$ block matrix 
$[M_{ij}]$, where $M_{ii} = 0$ for $1 \leq i \leq n$ and $M_{ij} =  2 a I + b D_{m}$ for 
$1 \leq i < j \leq n$. It is straightforward to verify that $E \in E(N, d)$ and Theorem 1 can be applied
to deduce $g \neq 0$.

To prove (ii), first note that $m n - r \geq 1$. Let $E \in E(N)$ be the $(n + 1) \times (n + 1)$ block matrix 
$[M_{ij}]$ defined as follows. For $1 \leq i \leq n + 1$, $M_{ii} = 0$. If $ m n - r \leq m$, then for 
$1 \leq i < j \leq n + 1$, $M_{1j}$ is the $(m n - r) \times m$ matrix having each entry equal to $c$ and 
$M_{ij} =  2 a I + b D_{m}$. If $m < m n - r$, then for $1 \leq i < j \leq n + 1$, $M_{ij} =  2 a I + b D_{m}$ 
and $M_{i (n + 1)}$ is the $m \times (m n - r)$ matrix having each entry equal to $c$. Then clearly $E \in E(N, d)$. 
If $m n - r = m $, then $m (m n - r) c = 2 m a + m (m -1) b$ is necessarily an even integer. Now it is 
straightforward to verify that Theorem 1 can be employed to infer $g \neq 0$.

To prove (iii), let $E \in E(N)$ be the $3 \times 3$ block matrix $[M_{ij}]$ such that $M_{rr} = 0$ for
$1 \leq r \leq 3$, $M_{12} = M_{21}^{T}$ is the $l \times m$ matrix having each entry equal to $a$, 
$M_{13} = M_{31}^{T}$ is the $l \times n$ matrix having each entry equal to $b$ and $M_{23} = M_{32}^{T}$ 
is the $m \times n$ matrix having each entry equal to $c$. By hypothesis, each of $a$, $b$, $c$ is a positive
integer. Since $d = m a + n b = l a + n c = l b + m c$, we have $E \in E(N, d)$. As before, it is easily 
verified that Theorem 1 is indeed applicable in this case and hence $g \neq 0$.

To prove (iv), let $m := 1$, $n := 4 u v + 2 v + 1$ and $r := 8 u v - 4 t v + 4 v$. Clearly, $n \geq 7$ and
$N = 2 m n - r$. Since $t \leq 2 u \leq 2 t - 1$, we have $1 \leq r \leq n - 1$. Define $a := (2 s + 1)(2 u - t + 1)$
and say $b := 1$. Then letting $c := (2 s + 1)(2 u + 1)$, we have $c \geq 3$ and $c r = (n - 1) [2 a + (m - 1) b]$.
Observe that the positive integers $a$, $b$, $c$, $m$, $n$, $r$ satisfy all the requirements of (ii). Thus, by taking
$E \in E(N, d)$ as described in the proof of (ii), we infer that $g \neq 0$. If $w$ denotes the weight of $g$, then
$2 w = N d$ and hence $w = (2 s + 1) (2 t v + 1) (2 u + 1) (4 u v + 2 v + 1)$. Since $w$ is an odd integer, $g$ is
a skew invariant. 

Lastly, to prove(v), suppose $E \in E(N, d)$ is such that each entry of $E$ is strictly less than $d$ and 
$Symm _{N} (\delta (z, E)) \neq 0$. Let $E^{*}$ be the $2 \times 2$ block matrix $[C_{ij}]$, where $C_{11} := 0$,
$C_{22} := E$ and $C_{12} = C_{21}^{T}$ is the $(N - 1) \times N$ matrix with each entry equal to $d$. Clearly,
$E^{*} \in E(2 N  - 1, d N)$ and Theorem 1 can be applied to infer $g \neq 0$. $\Box$ \\

\noindent \underline {\bf Example 2}: We continue assuming $N \geq 3$.
\begin{enumerate}
\item $N = 4 e$. Using (i) of Theorem 2 with $n := 2$ and $m :=  2 e$, we obtain nonzero invariants of degree 
$d$ for $d =  2 e + 1$ and all $d \geq N - 1$. If $char\, k = 0$ and $d \leq N - 2$ is even, then Lemma 2 yields 
a nonzero invariant of degree $d$. 
\item With the notation of (iii), let $Y := \{1 \leq d \in {\Bbb Z} \, \mid \, a, b, c \in {\Bbb Z} \}$ and 
\[ y \, := \, \frac{2 l m n}{gcd (N - 2 l,  N - 2 m, N - 2n, 2 l m n)} . \]
Then it is straightforward to verify that $d \in Y$ if and only if $d = s y$ for some positive integer $s$. 
Of course, $2 l m n \in Y$; but $y$ can be strictly less than $2 l m n$ ({\it e.g.}, consider $(l, m, n) := (2, 5, 6)$
or $(l, m ,n) := (9, 15, 21)$). If $l + m + n$ is odd and $d = 2 \; mod \; 4$, then the resulting $g$
is a nonzero skew invariant. So, (iii) produces skew invariants for binary forms of odd degrees (in contrast to (iv)). The 
least value of $N$ for which (iii) may be used to obtain skew invariants, is $N = 3 + 5 + 7 = 15$; whereas for the 
ones that can be obtained by using (iv) is $N = 2 (2 \cdot 2 \cdot 1 + 1) = 10$. For $3$-part partitions $N = l + m + n$
with $l \leq m \leq n < l + m$, by imposing additional requirements such as: $(l + m - n) d$ is divisible by $4$ if 
$l = m$ and so on, hypotheses of Theorem 1 can be satisfied. Assertion (iii) can be generalized for certain types 
of partitions of $N$ into $4$ or more parts; the task of formulating such generalizations is left to the reader.
\item Let $E \in \{ E_{1}, E_{2} \} \subset E(5, 18)$, where $E_{1}, E_{2}$ are as in the second example
above Theorem 2. For $2 \leq n \in {\Bbb Z}$, let $d_{n}$, $M_{n} \in E(2^{n} + 1, d_{n})$ be inductively 
defined by setting $d_{2} := 18$, $M_{2} := E$, $d_{n + 1} := (2^{n} + 1) d_{n}$ and where $M_{n+ 1} := M_{n}^{*}$,
is derived from $M_{n}$ as in (iv) of Theorem 2. Then by (v) of Theorem 2, 
$g_{n} := Symm _{2^{n} + 1} (\delta (z, M_{n}))$ is a nonzero skew invariant of a binary form of degree $2^{n} + 1$
for $2 \leq n \in {\Bbb Z}$.  \\
\end{enumerate}

\noindent \underline {\bf Remark 4}: Theorem 2 exhibits the simplest applications of Theorem 1. At present, there 
does not exist a characterization of pairs $(N, d)$ for which Theorem 1 can be used to obtain a nonzero invariant. 
Interestingly, it is impossible to use Theorem 1 to construct invariants corresponding to certain pairs $(N, d)$, 
{\it e.g}, consider $(N, d) = (5, 18)$: an elementary computation verifies that Hermite's invariant of a binary quintic 
can not be constructed via Theorem 1. A `good' generalization of Theorem 1, if it exists, should repair this failing. \\ 

\noindent \underline {\bf Definitions}: Let $n$, $s$ be a positive integers. 
\begin{enumerate}
\item Let $\preceq$ denote the lexicographic order on ${\Bbb Z}^{s+1}$.
\item For $\alpha := (a_{1}, \ldots , a_{s+1}) \in {\Bbb Z}^{s+1}$, let $|\alpha | := \sum _{i=1}^{s + 1} a_{i}$ and
\[wt(n, \alpha) \; := \; \frac{1}{2} \left [n^2 -  \left ( \sum _{i=1}^{s + 1} a_{i}^2 \right ) \right ]  . \]
\item Define $\wp(s, n) := (\wp_{1} (s, n), \ldots , \wp_{s + 1}(s, n)) \in {\Bbb Z}^{s + 1}$, where
\[ \wp_{j} (s, n) \; := \; \left \lfloor \frac{n - \sum_{1 \leq i \leq j - 1} \wp_{i}}{s + 2 - j} - \frac{(s + 1 - j)}{2} \right  \rfloor 
\;\;\;\;\;\; \mbox{for $1 \leq j \leq s + 1$.} \]
\item Let $\varpi (s, n) := wt(n, \wp(s, n))$. 
\item By $\Im (s, n)$ we denote the set of all $\alpha := (a_{1}, \ldots, a_{s + 1}) \in {\Bbb Z}^{s + 1}$
such that $a_{1} < a_{2} < \cdots < a_{s + 1}$ and $|\alpha| = n$. Let ${\Bbb P}(s, n)$ be the subset of $\Im (s, n)$
consisting of $(a_{1}, \ldots, a_{s + 1}) \in \Im (s, n)$ with $a_{1} \geq 1$. 
\item For $(i, j) \in {\Bbb Z}^2$ with $1 \leq i < j \leq s + 1$, let $\eta(i, j) := (\eta_{1}, \ldots, \eta_{s+1})$ 
where $\eta _{r} = 0$ if $r \neq i, j$, $\eta _{i} = 1$ and $\eta _{j} = -1$. An $(s + 1)$-tuple $\beta$ is said to be 
an {\it elementary modification} of $\alpha \in {\Bbb Z}^{s + 1}$ provided $\beta = \alpha + \eta (i, j)$ for some 
$1 \leq i < j \leq s + 1$. An $(s + 1)$-tuple $\beta$ is said to be a {\it modification} of $\alpha \in {\Bbb Z}^{s + 1}$ 
if there is a finite sequence $\alpha = \alpha_{1}, \ldots , \alpha _{r} = \beta$ such that $\alpha _{i}$ is an 
elementary modification of $\alpha _{i - 1}$ for $2 \leq i \leq r$.\\
\end{enumerate}

\noindent \underline {\bf Lemma 3}: Fix positive integers $n$, $s$ and let $e$ be the integer such that 
\[ n - \frac{s (s + 1)}{2} \; = \; \left \lfloor \frac{n}{s + 1} - \frac{s }{2} \right \rfloor (s + 1) + e . \]
Let $\wp(s, n) = (p_{1}, \ldots , p_{s + 1})$. Then, the following holds.
\begin{description}
\item[(i)] We have 
\[p _{j} \; = \; \left \{\begin{array}{ll} 
p _{1} + j - 1 & \mbox{if $1 \leq j \leq s + 1 - e$, and} \\
p _{1} + j & \mbox{if $s + 2 - e \leq j \leq s + 1$.} \end{array} \right . \]
In particular, $\wp (s, n) \in \Im (s, n)$. Moreover, if $(s + 1)(s + 2) \leq 2 n$, then 
$\wp (s, n) \in {\Bbb P}(s, n)$.
\item[(ii)] We have 
\begin{eqnarray*} 
\varpi (s, n) & = & 
\frac{(s + 1)(s + 2)}{2} \left \lfloor \frac{n}{s + 1} - \frac{s }{2} \right \rfloor ^{2} \\ & & +
\frac{(s + 1)^{2} (s + 2) - 2 n(s + 2)}{2} \left \lfloor \frac{n}{s + 1} - \frac{s }{2}  \right \rfloor \\  & & + 
\frac{3 (s + 1)^4 + 2 (s + 1)^3 - 3(1 + 4 n) (s + 1)^2 - 2 (1 + 6 n) (s + 1) + 24 n^2 }{24}. 
\end{eqnarray*}
\item[(iii)] Let  $\alpha := (a_{1}, \ldots, a_{s + 1}) \in \Im (s, n)$. Then, $\alpha \preceq \wp (s, n)$, 
$\wp (s, n)$ is a modification of $\alpha$ and 
\[ \sum _{1 \leq i < j \leq s + 1} a_{i} a_{j} \; = \; wt(n, \alpha) \; \leq \; \varpi (s, n). \]
\item[(iv)] ${\Bbb P}(s, n) \neq \emptyset$ if and only if 
\[ s  \;  \leq  \; \left \lfloor \frac{\sqrt{8 n + 1} - 1}{2} \right \rfloor - 1. \] 
\item[(v)] Suppose $s \geq 2$, $(s + 1)(s + 2) \leq 2 n$ and $p_{1} + e = b s + d$ where $b$, $d$ are nonnegative integers 
with $d \leq s - 1$. Then, letting $\wp (s - 1, n) := (q_{1}, \ldots , q_{s})$, we have $q_{1} = p_{1} + b + 1$ and
\[\varpi (s, n) - \varpi (s - 1, n) \; = \; p_{1} (s + 1 - e) + b d (s + 1) + \frac{1}{2} b (b - 1) s (s + 1). \]
In particular, $q_{1} > p_{1}$ and $\varpi (s, n) - \varpi (s - 1, n) \geq 2 p_{1}$. If $p_{1} = 1$, then $2 \leq q_{1} \leq 3$
and $2 \leq \varpi (s, n) - \varpi (s - 1, n) \leq s + 2$. 
\item[(vi)] Suppose $s \geq 2$, $(s + 1)(s + 2) \leq 2 n$ and let $v(s, n) := (v_{1}, \ldots , v_{s})$ where $v_{i} := i$
for $1 \leq i \leq s$ and $v_{s} = n - (1/2)s (s + 1)$. Then, $v(s, n) \preceq \alpha$ and $wt(n, v(s, n)) \leq wt(n, \alpha)$
for $\alpha \in {\Bbb P}(s, n)$. \\
\end{description}

\noindent \underline {\bf Proof}: Note that $0 \leq e \leq s$ and hence $s + 1 - e \geq 1$. Suppose 
$1 \leq j \leq s + 1 - e$ is such that $p_{i} = p_{1} + i - 1$ for $1 \leq i \leq j$. Then,
\begin{eqnarray*} 
p_{j + 1} & = & \left \lfloor p_{1} - \frac{j (j - 1) - s (s + 1)  - 2 e + (s - j)(s + 1 - j) }{2 (s + 1 - j)} \right  \rfloor \\
& = & \left \lfloor p_{1} + j + \frac{e}{s + 1 - j} \right  \rfloor. 
\end{eqnarray*}
If $j < s + 1 - e$, then $e < s + 1 - j$ and hence $p_{j+1} = p_{1} + j$. If $j = s + 1 - e$, then $p_{j+1} = p_{1} + j + 1$.
Next suppose (i) holds for some $j$ with $s + 2 - e \leq j \leq s$. Then,
\begin{eqnarray*} 
p_{j + 1} & = & 
\left \lfloor p_{1} - \frac{j (j - 1) - s (s + 1) + 2 (j + e - s - 1) - 2 e + (s - j)(s + 1 - j) }{2 (s + 1 - j)} \right \rfloor \\
& = & p_{1} + j + 1.
\end{eqnarray*}
Clearly, $p_{1} < p_{2} < \cdots < p_{s + 1}$ and if $(s + 1)(s + 2) \leq 2 n$, then 
$p_{1} \geq 1$. Also, $|\wp(s, n)| = p_{1} (s + 1) + [s (s + 1) / 2] + e = n$. Thus (i) holds. 

Let $u (X), v(X) \in {\Bbb Z}[X]$ be defined by 
\[ v(X) \; = \; \prod _{j = 0}^{s + 1} (X + p_{1} + j) \; = \; (X + p_{1} + s + 1 - e) u(X). \] 
Then, $\varpi (s, n)$ is the coefficient of $X^{s - 1}$ in $u(X)$. The coefficient of $X^{s}$ in $v(X - p_{1})$ is 
\[ \frac{1}{2} \left ( \sum _{i = 0}^{s + 1} i \right )^{2} -  \frac{1}{2} \sum _{i = 0}^{s + 1} i^2 \; = \; 
\frac{(3 s + 5) (s + 2) (s + 1) s}{24} \]. Now a straightforward computation verifies (ii).

Obviously, $wt(n, \alpha) < n^2$ for all $\alpha \in \Im(s, n)$. If $\beta \in \Im (s, n)$ is an elementary modification 
of $\alpha = (a_{1}, \ldots, a_{s + 1}) \in \Im(s, n)$, then note that $wt (n, \beta) > wt (n, \alpha)$. Hence $\alpha$
has a modification $v \in \Im (s, n)$ that is `final' in the sense that no member of $\Im (s, n)$ is an elementary modification 
of $v$. Fix such $v := (v_{1}, \ldots , v_{s+1})$. If $1 \leq i \leq s + 1$ is such that $v_{i+1} > v_{i} + 2$, then
$v + \eta (i, i + 1) \in \Im (s, n)$; this contradicts our assumption about $v$. So, $v_{i} + 1 \leq v_{i+1} \leq v_{i} + 2$ 
for all $1 \leq i \leq s$. If there are $1 \leq i < j \leq s + 1$ such that $v_{i+1} = v_{i} + 2$ as well as $v_{j+1} = v_{j} + 2$, 
then $v + \eta (i, j) \in \Im (s, n)$; an impossibility. Hence $a_{i + 1} = a_{i} + 2$ for at most one $i$ with 
$1 \leq i \leq s$. Consequently, $n = |v| = (s + 1) v_{1} + (s + 1 - j) + [s (s + 1) / 2]$ for some $j$ with
$1 \leq j \leq s + 1$. Clearly, $j = s + 1 - e$ and in view of (ii), we have $v = \wp(s, n)$. Thus $\wp (s, n)$ is a 
modification of $\alpha$. In particular, $wt(n, \alpha) \leq \varpi (s, n)$ and $\alpha \preceq \wp (s, n)$. The equality
displayed on the left in (iii) readily follows from the definition of $wt (n , \alpha)$. Thus (iii) holds. 

Assertion (iv) is simple to verify. To prove (v), assume $s \geq 2$ and let $p_{1} + e = b s + d$ where
$b$, $d$ are nonnegative integers with $d \leq s - 1$. Consequently, $q_{1} = p_{1} + b + 1 > p_{1}$. Using (ii) 
$\varpi(s, n) - \varpi (s - 1, n)$ can be computed in a straightforward manner. If $e \leq s - 1$, then 
$\varpi(s, n) - \varpi (s - 1, n)$ is clearly $\geq 2 p_{1}$. If $e = s$, then we have $b \geq 1$ and since 
$(b - 1) s = p_{1} - d$,
\[\varpi(s, n) - \varpi (s - 1, n) \; \geq \; p_{1} \left ( 1 + \frac {1}{2} b(s + 1) \right ) \; \geq \; 2 p_{1}. \]
If $p_{1} = 1$, then since $0 \leq e \leq s$ and $s \geq 2$, we have $0 \leq b \leq 1$. If $e \leq s - 2$, then $b = 0$
and hence $q_{1} = 2$, $\varpi(s, n) - \varpi (s - 1, n) = s + 1 - e \leq s + 1$. If $e =  s - 1$, then $b = 1$, $d = 0$
and hence $q_{1} = 3$, $\varpi(s, n) - \varpi (s - 1, n) = 2$. Lastly, if $e = s$, then $b = 1 = d$ and hence 
$q_{1} = 3$, $\varpi(s, n) - \varpi (s - 1, n) = s + 2$. This establishes (v). The proof of (vi) is left to 
the reader. $\Box$ \\

\noindent \underline {\bf Lemma 4}: Let $m, n , t \in {\Bbb Z}$ and $(b_{1}, \ldots , b_{m}) \in {\Bbb Z}^{m}$ be
such that $m \geq 1$, $n \geq 1$, $b_{1} + \cdots + b_{m} = t$ and $b_{i} \geq 0$ for $1 \leq i \leq m$. 
Let $t = q n + r$, where $q$, $r$ are integers with $q \geq 0$ and $0 \leq r < n$. Then, there exists an $m \times n$ 
matrix $A := [a_{ij}]$ satisfying the following.
\begin{description} 
\item[(i)] $0 \leq a_{ij} \in {\Bbb Z}$ for $1 \leq i \leq m$, $1 \leq j \leq n$ and $\|A \| = t$.
\item[(ii)] \[c_{j}(A) \; := \; r_{j}\left (A^{T} \right ) \; = \; \left \{\begin{array}{ll} 
q + 1 & \mbox{if $1 \leq j \leq r$ and} \\ q & \mbox{if $r + 1 \leq j \leq n$.} \end{array} \right . \]
\item[(iii)] $r_{i}(A) = b_{i}$ for $1 \leq i \leq m$. \\
\end{description}

\noindent \underline {\bf Proof}: Let $t = q n + r$, where $q$, $r$ are integers with 
$q \geq 0$ and $0 \leq r < n$. Our proof proceeds by induction on $m$. If $m = 1$, then let $a_{1j} := q + 1$ if 
$1 \leq j \leq r$ and $a_{1j} := q$ if $r + 1 \leq j \leq n$. Henceforth suppose $m \geq 2$ and $b_{m} = \ell n + \rho$ 
where $\ell$, $\rho$ are integers with $\ell \geq 0$ and $0 \leq \rho < n$. 

\noindent Case 1: $\rho \leq r$. By our induction hypothesis there is an $(m - 1) \times n$ matrix $[a_{ij}]$ 
such that $0 \leq a_{ij} \in {\Bbb Z}$ for $1 \leq i \leq m - 1$ and $1 \leq j \leq n$, $\|A \| = t - b_{m}$,
$a_{1j} + \cdots + a_{(m-1)j} = q - \ell + 1$ for $1 \leq j \leq r - \rho$, $a_{1j} + \cdots + a_{(m-1)j} = q - \ell$ 
for $r - \rho + 1 \leq j \leq n$ and $a_{i1} + \cdots + a_{in} = b_{i}$ for $1 \leq i \leq m - 1$. Define 
$a_{mj} := \ell$ for $1 \leq j \leq r - \rho$, $a_{mj} := \ell + 1$ for $r - \rho + 1 \leq j \leq r$ and 
$a_{mj} := \ell$ for $r + 1 \leq j \leq n$. Then, the resulting $m \times n$ matrix $[a_{ij}]$ is clearly the
desired matrix $A$.

\noindent Case 2: $\rho > r$. At the outset observe that $r < n + r - \rho < n$. As before, our induction hypothesis 
assures the existence of an $(m - 1) \times n$ matrix $[a_{ij}]$ such that $0 \leq a_{ij} \in {\Bbb Z}$ for 
$1 \leq i \leq m - 1$ and $1 \leq j \leq n$, $\|A \| = t - b_{m}$, $a_{1j} + \cdots + a_{(m-1)j} = q - \ell$ for 
$1 \leq j \leq n + r - \rho$, $a_{1j} + \cdots + a_{(m-1)j} = q - \ell - 1$ for $n + r - \rho + 1 \leq j \leq n$ and 
$a_{i1} + \cdots + a_{in} = b_{i}$ for $1 \leq i \leq m - 1$. Define $a_{mj} := \ell + 1$ for $1 \leq j \leq r$, 
$a_{mj} := \ell$ for $r + 1 \leq j \leq n + r - \rho$ and $a_{mj} := \ell + 1$ for $n + r - \rho + 1 \leq j \leq n$.
Then, the resulting $m \times n$ matrix $[a_{ij}]$ is the desired matrix $A$. $\Box$ \\

\noindent \underline {\bf Definitions}: Let $n$ and $w$ be positive integers. 
\begin{enumerate}
\item Define 
\[ \beta (n) \; := \; \left  \lfloor  \frac{ \sqrt{8 n + 1} - 1}{2} \right  \rfloor . \]
\item For an integer $s$ with $1 \leq s \leq \beta (n) - 1$ and an 
${\frak a} := (m_{1}, \ldots, m_{s+1}) \in {\Bbb P}(s, n)$, define
\[ \nu (w, {\frak a}) \; := \; {s - 1 + w - wt(n, {\frak a}) \choose s - 1} \]
and
\[ d (w, {\frak a}) \; := \; \left \{\begin{array}{ll} 
n - 1 + w - wt(n, {\frak a}) & \mbox{if $m_{1} = 1$,} \vspace{0.1in}\\
n - 1 + w - wt(n, {\frak a}) & \mbox{if $w = 1 + wt(n, {\frak a})$,} \vspace{0.1in}\\
n - m_{1} + 1 + \left \lceil \frac{w - wt(n, {\frak a})}{m_{1}} \right \rceil & \mbox{otherwise.} 
\end{array} \right . \]
\item Let $\nu (w, s, n) := \nu (w, \wp (s, n))$ and $d (w, s, n) := d (w, \wp (s, n))$. \\
\end{enumerate}

\noindent \underline {\bf Theorem 3}: Assume that $N$ is an integer $\geq 3$ and $k$ is a field of 
characteristic either $0$ or strictly greater than $N$. Let $F$ be the generic binary form of degree $N$ 
(as in the introduction). Let $s$ be an ineteger with $1 \leq s \leq \beta (N) - 1$ and let 
${\frak a}:= (m_{1}, \ldots, m_{s+1}) \in {\Bbb P}(s, N)$. Let $m := m_{1}$ and let $w$ be an integer such that 
$\theta := w - wt (N, {\frak a}) \geq 1$. Then, for a positive integer $d \geq d (w, {\frak a})$, there exist 
$\nu (w, {\frak a})$ $k$-linearly independent semi-invariants of $F$ of weight $w$ and degree $d$. \\

\noindent \underline {\bf Proof}: Fix an ordered $s$-tuple 
$(\theta_{1}, \ldots , \theta_{s})$ of nonnegative integers with 
\[ \theta_{1} + \cdots + \theta_{s} \; =  \; \theta .  \] Since 
$\theta \geq 1$, using Lemma 4 we obtain an $s \times m$ matrix $B^{*} := [b^{*}_{ij}]$ having nonnegative integer entries 
such that $r_{i} (B^{*}) = \theta_{i}$ for $1 \leq i \leq s$ and 
\[ \lfloor \theta / m \rfloor \; \leq \; c_{m} (B^{*}) \; \leq \; \cdots \; \leq \; c_{1} (B^{*}) \; =  \;
\lceil \theta / m \rceil. \] Let $u$ be the greatest positive integer such that $c_{u}(B^{*}) \geq 1$ and let $v$ be the 
least positive integer with $b^{*}_{vu} \geq 1$. Define an $s \times m$ matrix $B := [b_{ij}]$ as follows. 
If $u = 1$ (in particular, if $m = 1$), let $B = B^{*}$. If $u \geq 2$, then let $b_{ij} := b^{*}_{ij}$ for
$(i, j) \neq (v, 1), \, (v, u)$, let $b_{vu} := b^{*}_{vu} - 1$ and let $b_{v1} := b^{*}_{v1} + 1$. Then, $B$ has 
nonnegative integer entries, $r_{i} (B) = \theta_{i}$ for $1 \leq i \leq s$,  
\[\begin{array}{l} c_{1}(B) \; = \; \min \left \{ 1 + \lceil \theta / m \rceil, \, \theta \right \}, \;\;\;\mbox{and}\\
\lfloor \theta / m \rfloor  - 1 \; \leq \; c_{j}(B) \; \leq \;  \lceil \theta / m \rceil, \;\;\;
\mbox{for $2 \leq j \leq m$.} \end{array} \]
Using Lemma 4 again, we obtain matrices $A_{1}, \ldots , A_{s}$ with nonnegative integer entries such that 
\[ \begin{array}{ll} (1) & \mbox{$A_{l}$ has size $m \times m_{l + 1}$ for $1 \leq l \leq s$,} \\
(2) & \mbox{$r_{i} (A_{l}) = b_{li}$ for $1 \leq l \leq s$, $1 \leq i \leq m$ and} \\
(3) & \mbox{$\lfloor \theta_{l} / m \rfloor \leq c_{j} (A_{l}) \leq c_{j - 1} (A_{l}) \leq  \lceil \theta_{l} / m \rceil $
for $2 \leq j \leq m_{l + 1}$.}  \end{array} \]
Clearly, $\| A_{l} \| = \theta_{l}$ for $1 \leq l \leq s$. Furthermore, we have 
\[\begin{array}{ll} (4) & \mbox{$r_{1}(A_{1}) + \cdots + r_{1}(A_{s}) =  
\min \left \{ 1 + \lceil \theta / m \rceil, \, \theta \right \}$, and} \\
(5) & \mbox{ $r_{i}(A_{1}) + \cdots + r_{i}(A_{s}) \leq \lceil \theta / m \rceil$ for $2 \leq i \leq m$.}
\end{array} \]
Let ${\Bbb I}$ denote a matrix (of any chosen size) having each entry $1$. Let $M := [M_{ij}]$ be an $(s + 1) \times (s + 1)$ 
block-matrix such that $M_{ji}$ is the transpose of $M_{ij}$ for $1 \leq i \leq j \leq s + 1$, and the block $M_{ij}$ is 
a $m_{i} \times m_{j}$ matrix defined by 
\[M_{ij} \; := \; \left \{\begin{array}{ll} 0 & \mbox{$i = j$,} \\ 
{\Bbb I} + A_{j - 1} & \mbox{if $i = 1 < j \leq s + 1$,} \\
{\Bbb I} & \mbox{ if $ 2 \leq i < j \leq s + 1$. } \end{array} \right . \]
Let $M^{\prime}$ denote the $(N - 1) \times (N - 1)$ matrix obtained from $M$ by deleting the first row as well as the first
column of $M$. Then, $M \in E(N)$ and $M^{\prime} \in E(N - 1)$. Also, in view of properties (1) - (5), it is straightforward 
to verify that 
\[r_{1} (M) \; = \; d (w, {\frak a}) \; > \; r_{i} (M) \;\;\;\;\; \mbox{for $2 \leq i \leq N$,} \]
and each of $M$, $M^{\prime}$ satisfies requirements (1), (2), (i) - (iv) of Theorem 1. Hence letting 
$\phi (\theta_{1}, \ldots , \theta_{s}) := Symm_{N} \left (\delta (z, M) \right )$, we have 
$\phi (\theta_{1}, \ldots , \theta_{s}) \; \neq \; 0$ as well as $Symm_{N - 1} \left (  \delta (z, M^{\prime}) \right ) \; \neq \; 0$.
Observe that the coefficient of $z_{1}^{d (w, \alpha)}$ in $\phi (\theta_{1}, \ldots , \theta_{s})$ is the symmetrization of
$\delta (z^{\prime}, M^{\prime})$ where $z^{\prime} := (z_{2}, \ldots , z_{N})$. Since 
$Symm_{N - 1} \left (\delta (z, M^{\prime}) \right ) \neq 0$, we conclude that the $z_{1}$-degree (and hence also each $z_{i}$-degree) 
of $\phi (\theta_{1}, \ldots , \theta_{s})$ is exactly $d (w, {\frak a})$. Let $\alpha$ be the $k$-monomorphism employed in
Theorem 1. Then, as noted in no. 2 of Remarks 3, the $t$-initial coefficient of 
$\alpha \left ( \phi (\theta_{1}, \ldots , \theta_{s}) \right )$ is a nonzero constant ({\it i.e.}, element of $k$) multiple of
\[\eta (\theta_{1}, \ldots , \theta_{s}) \; := \; \prod _{1 \leq j \leq s} (t_{1} - t_{j + 1})^{\theta_{j}}
\prod _{1 \leq i < j \leq s + 1} (t_{i} - t_{j})^{m_{i} m_{j}}. \]
The set of all $\eta (\theta _{1}, \ldots , \theta _{s})$ ranging over the allowed choices of $s$-tuples 
$(\theta _{1}, \ldots , \theta _{s})$, is clearly a $k$-linearly independent subset of $k[t_{1}, \ldots , t_{s+1}]$.
Hence the corresponding set $S (\theta) $ of $\phi (\theta_{1}, \ldots , \theta_{s})$ is also a $k$-linearly independent subset 
of $k[z_{1}, \ldots , z_{N}]$. Of course $S (\theta) \subset k [y_{1}, \ldots , y_{N - 1}] \subset k[e_{1}, \ldots , e_{N}]$ 
(where $y_{1}, \ldots , y_{N -1}$ and $e_{1}, \ldots , e_{N}$ are as in the introduction). Given $\phi \in S (\theta)$, we 
homogenize $\phi$ to get a homogeneous polynomial of degree $d (w, {\frak a})$ in $a_{0}, \ldots , a_{N}$ as in the introduction. In this
manner we obtain a $k$-linearly independent set ${\Bbb S} (\theta)$ of semi-invariants of $F$ of degree $d (w, {\frak a})$ and weight $w$.
Obviously, $|{\Bbb S} (\theta) | = | S (\theta) | = \nu (w, {\frak a})$. Letting $ v := d - d (w, {\frak a})$, it follows that the set
$\left \{ a_{0}^{v} \sigma \; \middle | \; \sigma \in  {\Bbb S} (\theta) \right \}$ is also $k$-linearly independent. $\Box$ \\

\noindent \underline{\bf Example 3}: Here we consider the case of $3 \leq N \leq 7$. It is essential to point out that the lower 
bounds proved in [4], [12], [19] assume $N \geq 8$. To the best of our knowledge, there is nothing in the existing literature with
which we can compare the bounds in examples below. 
\begin{enumerate}
\item If $N = 3$, then $s = 1$ and $\varpi (1, 3) = 2$. In this case, Theorem 3 implies that for $0 \leq n \in {\Bbb Z}$, there 
exists a nonzero semi-invariant (of a binary cubic form $F$) of weight $2 + n$ and degree at least $2 + n$.
\item If $N = 4$, then $s = 1$ and $\varpi (1, 4) = 3$. In this case, Theorem 3 implies that for $0 \leq n \in {\Bbb Z}$, there 
exists a nonzero semi-invariant (of a binary quartic form $F$) of weight $3 + n$ and degree at least $3 + n$.
\item If $N = 5$, then $s = 1$ and $\varpi (1, 5) = 6$. In this case, Theorem 3 implies that for $0 \leq n \in {\Bbb Z}$, there 
exists a nonzero semi-invariant (of a binary quintic form $F$) of weight $6 + n$ and degree at least $4 + \lceil n / 2 \rceil$. Note 
that for the partition $1 < 4$, we can use Theorem 1 to verify the existence of a nonzero semi-invariant of weight $4 + n$ and degree 
at least $4 + n$. So, we obtain two $k$-linearly independent semi-invariants of weight $6 + n$ and degree at least $6 + n$. 
\item Assume $N = 6$. Then $1 \leq s \leq 2$, $\varpi (1, 6) = 8$ and $\varpi (2, 6) = 11$. Taking $s = 1$ in Theorem 3, we
infer the existence of a nonzero semi-invariant (of a binary sextic form $F$) of weight $8 + n$ and degree at least $8 + n$ for all
$0 \leq n \in {\Bbb Z}$. Next, taking $s = 2$, Theorem 3 assures the existence of $5 + n$ $k$-linearly independent semi-invariants of 
weight $16 + n$ and degree at least $10 + n$ for all $0 \leq n \in {\Bbb Z}$. 
\item Assume $N = 7$. Then $1 \leq s \leq 2$, $\varpi (1, 7) = 12$ and $\varpi (2, 7) = 14$. Letting $s = 1$ in Theorem 3,
we obtain a nonzero semi-invariant (of a binary heptic form $F$) of weight $12 + n$ and degree at least $5 + \lceil n / 3 \rceil$ 
for $0 \leq n \in {\Bbb Z}$. Using Theorem 1 for the partition $2 < 5$, we infer the existence of a nonzero semi-invariant of weight 
$10 + n$ and degree at least $6 + \lceil n / 2 \rceil$ for all $0 \leq n \in {\Bbb Z}$. Letting $s = 2$ in Theorem 3, we deduce the 
existence of $5 + n$ $k$-linearly independent semi-invariants of weight $18 + n$ and degree at least $5 + \lceil (n + 4) / 3 \rceil$ 
for all $0 \leq n \in {\Bbb Z}$. \\
\end{enumerate}

\noindent \underline{\bf Remark 5}: Let $N$, $w$ and $d$ are positive integers. Let 
\[PP(N, w, d) \; := \; \left \lceil \frac{4}{1000} \cdot (\min \{2 w,\, d^{2}, \, N^{2} \})^{\frac{-9}{4}} \cdot 
2^{\sqrt{\min \{2 w,\, d^{2} \, N^{2} \}}} \right \rceil .\]
If $\min \{N, d \} \geq 8$ and $w \leq N d / 2$, then by Theorem 1.2 of [12], there are at least $PP(N, w, d)$ $k$-linearly
independent semi-invariants (of a binary $N$-ic form $F$) of degree $d$ and weight $w$. Observe that for $(w, d)$ with 
$w  \geq N^{2} / 2$ and $d \geq N$, the bound $PP(N, w, d)$ is independent of $(w, d)$ ({\it i.e.}, depends only on $N$). 
In contrast, the lower bound $\nu (w, {\frak a})$ is a polynomial of degree $s - 1$ in $w$. The reader may wish to make 
similar comparison with results of [4].\\

\noindent \underline{\bf Example 4}: Let $\nu (w, N)  := \nu (w, \beta (N) - 1, N)$. Consider the case of $N = 15$.
Note that $\beta (N)  = 5$ and ${\Bbb P}(4, 15) = \{ \wp (4, 15) \}$. We have $\varpi (4, 15) = 85$ 
and $\wp _{1} (4, 15) = 1$. Let $\nu (w)  := \nu (w, 4, 15)$. Then, Theorem 3 assures that for $0 \leq n \in {\Bbb Z}$, 
we have at least $\nu  (85 + n)$ $k$-linearly independent semi-invariants of weight $85 + n$ and degree $ d \geq 14 + n$. 
Observe that $2 (85 + n) < (14 + n)^{2}$ for $n \geq 0$, $N^{2} = 225 < 2 (85 + n)$ for $n \geq 28$ and
\[\nu (85 + n) \; = \; {3 + n  \choose 3} \; = \; \frac{1}{6} n^3 + n^2 + \frac{11}{6} n + 1 \;\;\;\mbox{for $n \geq 0$.} \]
A straightforward computation verifies that $PP(15, 85 + n, d) = 1 < \nu (85 + n) $ for all  $n \geq 0$ and $d \geq 14 + n$. 
Let $semdim(w, d, N)$ denote the dimension of the $k$-vector space of semi-invariants (of our $N$-ic form $F$) 
of weight $w$ and degree $d$. Assume $k$ has characteristic $0$. Then, in the notation of the introduction, $semdim(w, d, N)$ is 
\[p_{w} (N, d) - p_{w-1} (N, d) \; := \; 
\mbox{the coefficient of $q^w$ in} \;\;(q - 1) {N + d \choose d}_{q}.  \]
The table below presents a MAPLE computation of  $\nu (85+n)$ and  $semdim(85 + n, 14 + n, 15)$ (denoted by $semdim$) 
for a small sample of values of the weight $w$ ({\it i.e.}, values of  $n$).
\begin{center}
\begin{tabular}{|c|c|c|} 
\hline 
$w$ & $\nu (w)$ & $semdim$ \\
\hline
& & \\
95 & 286 & 1020697 \\
& & \\
105 & 1771 & 4232793 \\
& & \\
115 & 5456 & 11374824 \\
\hline 
\end{tabular}
\quad
\begin{tabular}{|c|c|c|} 
\hline 
$w$ & $\nu (w) $ & $semdim$ \\
\hline
& & \\
125& 12341 & 25995316 \\
& & \\
135 & 23426 & 54621331 \\
& & \\
145 & 39711 & 108639772 \\
\hline 
\end{tabular}
\end{center}

Let $s = 3$ and ${\frak a} := v(3, 15) = (1, 2, 3, 9)$. Then, for integers $n \geq 0$, we have 
$\nu(65 + n, {\frak a}) = (1 / 2) (n + 2)(n + 1)$ and $d(65 + n, {\frak a}) = 14 + n$. At the other
extreme, if ${\frak a} = \wp (3, 15)$, then $\varpi (3, 15) = 80$ and $\wp _{1} (3, 15) = 2$. So, 
$\nu(80 + n, 3, 15) = (1 / 2) (n + 2)(n + 1)$ and $d(80 + n, 3, 15) = 14 + \lceil n / 2 \rceil$ for 
all $n \geq 0$. Thus for weights $65 \leq w < 80$, our lower bound is for degrees $\geq w - 1$; whereas, 
for weights $w \geq 80$ our lower bound is for degrees $\geq 14 + \lceil (w - 80) / 2 \rceil$. 
If $s = 2$, then $\varpi (2, 15) = 74$ and $\wp _{1} (2, 15) = 4$. Hence $\nu(74 + n, 2, 15) = n + 1$ and 
$d(74 + n, 2, 15) = 12 + \lceil n / 4 \rceil$ for all $n \geq 0$. For $s = 1$, we have $\varpi (1, 15) = 56$ and 
$\wp _{1} (1, 15) = 7$. Consequently, $\nu(56 + n, 1, 15) = 1$ and $d(56+ n, 1, 15) = 9 + \lceil n / 7 \rceil$.\\ \\

\noindent{\bf References} 
\begin{description}
\item[1.] Abhyankar, Shreeram. S. (1988) {\em Enumerative combinatorics of Young tableaux}, Monographs and Textbooks 
in Pure and Applied Mathematics, 115. Marcel Dekker, Inc., New York.
\item[2.] Benjamin, Arthur T.; Quinn, Jennifer J.; Quinn, John J.; W{\'o}js, Arkadiusz (2001){\em Composite fermions and 
integer partitions}, J. Combin. Theory Ser. A 95, no. 2, 390 -397.
\item[3.] Brennan, Joseph P. (1997) {\em Invariant theory in characteristic p: Hazlett's symbolic method for binary quantics},
Factorization in integral domains (Iowa City, IA, 1996), 257-269, Lecture Notes in Pure and Appl. Math., 189, Dekker, New York.
\item[4.] Dhand, Vivek {\em A combinatorial proof of strict unimodality for q-binomial coefficients}, Discrete 
Math. 335 (2014), 20-24.
\item[5.] Dunajski, M.; Penrose, R. (2017) {\em On the quadratic invariant of binary sextics}, Math. Proc. 
Cambridge Philos. Soc., 162, no. 3, 435-445.
\item[6.] Elliot, E. B. (1913) {\em An Introduction to the Algebra of Quantics}, Chelsea Publishing
Company, New York, 1964, Second edition, reprint.
\item[7.] Grace J. H.; Young, A. (1903) {\em The Algebra of Invariants}, Chelsea Publishing Company, 
New York, 1964, reprint.
\item[8.] Kung, J. P. S. and Rota, G.-C. (1984) {\em The invariant theory of binary forms}, Bull. Amer. Math. Soc., 10: 
27 - 85.
\item[9.] Mulay, S. (2018) {\em Graph-monomials and invariants of binary forms}, arXiv:1809.00369.
\item[10.] Mulay, S.; Quinn, John J.; Shattuck, M. (2016) {\em Correlation diagrams: an intuitive approach
to correlations in quantum Hall systems}, Journal of Physics: Series C, 702, (012007 - 1) - (012007 - 9). 
\item[11.] Mulay, S.; Quinn, John J.; Shattuck, M. (2018) {\em Strong Fermion Interactions in Fractional Quantum
Hall States, Correlation Functions}, Springer Series in Solid-State Sciences, 193. 
\item[12.] Pak, Igor; Panova, Greta (2017) {\em Bounds on certain classes of Kronecker and q-binomial coefficients},
J. Combin. Theory Ser. A 147, 1-17.
\item[13.] Petersen, J. (1891) {\em Die Theorie der regularen Graphs}, Acta Math., 15, 193-220.
\item[14.] Petersen, J. (1897) {\em Theorie des equations Algtbriques}, Gauthier-Villars, Paris.
\item[15.] Quinn, John J.; W\'ojs, A. (2000) {\em Composite fermions in fractional quantum Hall systems}, 
Journal of Physics: Condensed Matter 12, R265-R298.
\item[16.] Sabidussi, G. (1992) {\em Binary invariants and orientations of graphs}, Discrete Math., 101, 251-277.
\item[17.] Sylvester, J. J. (1878) {\em On an application of the New Atomic Theory to the graphical representation of
invariants and covariants of binary quantics, with three appendices}, Amer. J. Math. 1, 64-125.
\item[18.] Sylvester, J. J. (1878) {\em Proof of the hitherto undemonstrated fundamental theorem of invariants}, 
Phi-los. Mag. 5, 178-188 (reprinted in: Coll. Math. Papers, vol. 3, Chelsea, New York, 1973). 
\item[19.] Zanello, Fabrizio (2015) {\em Zeilberger's KOH theorem and the strict unimodality of q-binomial coefficients}, 
Proc. Amer. Math. Soc. 143, no. 7, 2795-2799.
\end{description}

\end{document}